\documentstyle[12pt]{article}

\let \ttorg \tt \def \tt{\ttorg \obeyspaces}

\begin{document}

\date{}

\title{\Large\bf Oriented Quantum Algebras, Categories and Invariants of Knots and Links}
\author{Louis Kauffman\thanks{Research supported in part by NSF Grant
DMS 920-5227} $\;\;$
and $\;\;$
David E. Radford\thanks{Research supported in part by NSF Grant
DMS 980 2178}\\ 
Department of Mathematics, Statistics \\
and Computer Science (m/c 249)    \\
851 South Morgan Street   \\
University of Illinois at Chicago\\
Chicago, Illinois 60607-7045}

%\author{Louis H. Kauffman  and David E. Radford\\ Department of
%Mathematics, Statistics and Computer Science \\ University of
%Illinois at Chicago \\ 851 South Morgan Street\\ Chicago, IL,
%60607-7045}

\maketitle

\thispagestyle{empty}

\subsection*{\centering Abstract}

{\em This paper defines the concept of an oriented quantum algebra and 
develops its application to the construction of quantum link invariants.
We show, in fact, that all known quantum link invariants can be put into this framework.}

\section{Introduction} An oriented quantum algebra $(A, \rho, D, U)$ is an abstract
model for an oriented quantum link invariant. This model is based
on a solution to the Yang-Baxter equation and some extra structure
that serves to make an invariant possible to construct. The
definition of an oriented quantum algebra is as follows:   We are
given an algebra $A$ over a base ring $k$, an invertible solution
$\rho$ in $A \otimes A$ of the Yang-Baxter equation (in the
algebra formulation of this equation -- see the Remark below), and automorphisms $U:A
\longrightarrow A$ and $D:A \longrightarrow A$ of the algebra. It
is assumed that $D$ and $U$ commute and that

$$(U \otimes U)\rho = \rho$$

\noindent and $$(D \otimes D)\rho = \rho$$

\noindent and that 
$$[(1_{A} \otimes U)\rho)][(D \otimes 1_{A^{op}})\rho^{-1}] 
= 1_{A \otimes A^{op}}.$$ 

\noindent and 

$$[(D \otimes 1_{A^{op}})\rho^{-1}][(1_{A} \otimes U)\rho)]
= 1_{A \otimes A^{op}}.$$ 

\noindent In other words, $[(1_{A} \otimes U)\rho)]$ and $[(D \otimes 1_{A^{op}})\rho^{-1}]$
are inverses in the algebra $A \otimes A^{op}.$  
\vspace{3mm}

\noindent Here $A^{op}$ denotes the opposite algebra. The first
equation is formulated in the tensor product of $A$ with itself,
while the second equation is formulated in the tensor product of
$A$ with its opposite algebra. \vspace{3mm}

When $U=D=T$, then $A$ is said to be {\em balanced}.  In this case
$$(T \otimes T)\rho = \rho,$$

$$[(1_{A} \otimes T)\rho)][(T \otimes 1_{A^{op}})\rho^{-1}] = 1_{A
\otimes A^{op}}$$

\noindent and 

$$[(T \otimes 1_{A^{op}})\rho^{-1}][(1_{A} \otimes T)\rho)] = 1_{A
\otimes A^{op}}.$$

In the case where $D$ is the identity mapping, we call the 
oriented quantum algebra {\em standard}.  As we shall see in section 6, 
the invariants defined by Reshetikhin and Turaev (associated with a 
quasi-triangular Hopf algebra) arise from standard oriented quantum algebras.
It is an interesting structural feature of algebras that we have elsewhere 
\cite{GAUSS} called {\em quantum algebras} (generalizations of quasi-triangular Hopf algebras)
that they give rise to standard oriented quantum algebras.
\vspace{3mm}

We shall see that appropriate matrix representations of oriented quantum 
algebras or the existence of certain traces on these algebras
allow the construction of oriented invariants of knots and links.
These invariants include all the known quantum link
invariants at the time of this writing. \vspace{3mm}

\noindent {\bf Remark.} Note that we have the Yang-Baxter elements $\rho$ and $\rho^{-1}$
in $A \otimes A.$ We assume that $\rho$ and $\rho^{-1}$ satisfy the algebraic Yang-Baxter
equation. This equation (for $\rho$) states

$$\rho_{12}\rho_{13}\rho_{23} = \rho_{23}\rho_{13}\rho_{12}$$

\noindent where $\rho_{ij}$ denotes the placement of the tensor factors of $\rho$ in the 
$i$-th and $j$-th tensor factors of the triple tensor product $A \otimes A \otimes A.$
\vspace{3mm}

\noindent We write $\rho = \Sigma e \otimes e'$ and
$\rho^{-1} = \Sigma E \otimes E'$ to indicate that these elements are sums of
tensor products of elements of $A$. The expression $e \otimes e'$ is thus a generic 
element of the tensor product. However, we often abbreviate and write
$\rho = e \otimes e'$ and $\rho^{-1} = E \otimes E'$ where the summation is implicit. 
We refer to $e$ and $e'$ as the {\em signifiers} of $\rho$,
and $E$ and $E'$ as the {\em signifiers} of $\rho^{-1}$. For example, 
$\rho_{13} = e \otimes 1 \otimes e'$ in $A \otimes A \otimes A.$ 
\vspace{3mm}

\noindent Braiding operators, as they appear in knot theory, differ from the algebraic 
Yang-Baxter elements by a permutation of tensor factors. This point is crucial to the 
relationship of oriented quantum algebras and invariants of knots and links, and will 
be dealt with in the body of the paper (specifically in Section 3).
\vspace{3mm} 

\noindent {\bf Remark.} In writing a function $F(x)$ we often abbreviate this 
expression to $Fx$ when this entails no ambiguity.  Thus $F(G(H(J(x)))) = FGHJx$
in this notation.
\vspace{3mm}

\noindent {\bf Remark.} If $A$ is an  $n \times n$ matrix and $v$ is a row vector of 
dimension $n$, then $vA$ is a row vector of dimension $n$. If matrices are regarded 
as functions on the left ( $(v)A = vA$), then the order of composition proceeds from 
right to left as in $((v)A)B = vAB.$ In some of our applications it is useful to think of
the matrices as acting on the left in this fashion so that $AB$ denotes both matrix 
product and composition of functions. The role of this remark will become clear in the
sections on matrix models, state models and quasi-triangular Hopf algebras.
\vspace{3mm} 

The paper is organized as follows. Section 2 describes the tangle category $Tang$ and the 
flat tangle category $Flat$. Section 3 describes the category $Cat(A)$
of an oriented quantum algebra $A$ and the functor $F: Tang \longrightarrow Cat(A).$
In fact, this section motivates the definition of an oriented quantum algebra by 
showing just what algebraic conditions are neccessary for the functor $F$ to be an
invariant of regular isotopy of tangles. This analysis is how we discovered the definition 
of these algebras. Section 4 discusses matrix models for the invariants of knots and links
associated with oriented quantum algebras. Included in Section 4 is a discussion of bead
sliding for reformulating the evaluations of the matrix models and a specific example of a
matrix oriented quantum algebra that gives rise to specializations of the Homfly polynomial.
Section 5 is a discussion of combinatorial state sum invariants of knots and links. The 
contents of this section can be used to verify the relationship of the algebra of the 
previous section with the Homfly polynomial. Section 6 shows that representations of 
quasitriangular ribbon Hopf algebras have the stucture of oriented quantum algebras.
We use this section to show explicitly how the link invariants of Reshtikhin and Turaev
fit into our framework.
\vspace{3mm}

\section{The Tangle Category} We recall the {\em oriented tangle
category} denoted here by $Tang$.  This is a category that
formalizes the structure of knot and link diagrams in a manner
suitable for the construction of quantum link invariants. (The reader 
should note that $Tang$ refers, as explained below, to tangles will all
multiplicities of input and output.)
\vspace{3mm}

We shall refer to embeddings of disjoint unions of circles and
arcs into three dimensional space as {\em string}. In this paper,
all strings will be {\em oriented}. This means that each string is
equipped with a preferred direction, usually indicated in a
diagram by an arrowhead drawn on the string. When we speak of
matching strings in tangles to compose them (see below and Figure 0), we assume
that the strings have compatible orientations so that the
composition is also oriented. \vspace{3mm}

 {\tt    \setlength{\unitlength}{0.92pt}
\begin{picture}(391,206)
\thinlines    \put(87,163){\line(-3,2){12}}
              \put(59,148){\line(2,1){27}}
              \put(42,137){\line(5,3){11}}
              \put(33,149){\line(3,-4){9}}
              \put(54,165){\line(-4,-3){21}}
              \put(27,192){\line(3,-5){16}}
              \put(49,157){\line(3,-5){10}}
              \put(59,140){\line(-3,-2){18}}
              \put(41,128){\line(-1,1){18}}
              \put(23,148){\line(3,5){12}}
              \put(43,175){\line(1,0){18}}
              \put(55,164){\line(3,-2){11}}
              \put(70,143){\line(0,-1){27}}
              \put(65,167){\line(6,-5){9}}
              \put(79,150){\line(1,-2){10}}
              \put(89,129){\line(0,-1){28}}
              \put(71,117){\line(-1,0){35}}
              \put(36,116){\line(0,-1){19}}
              \put(65,167){\line(3,5){14}}
              \put(10,107){\framebox(99,74){}}
              \put(149,128){\framebox(57,55){}}
              \put(235,126){\framebox(57,55){}}
              \put(163,183){\line(0,1){13}}
              \put(183,183){\line(0,1){13}}
              \put(162,130){\line(0,-1){14}}
              \put(184,129){\line(0,-1){14}}
              \put(272,126){\line(0,-1){14}}
              \put(248,126){\line(0,-1){14}}
              \put(274,180){\line(0,1){13}}
              \put(250,181){\line(0,1){13}}
              \put(165,144){\makebox(25,23){A}}
              \put(249,142){\makebox(17,24){B}}
              \put(338,76){\makebox(17,24){B}}
              \put(339,143){\makebox(25,23){A}}
              \put(339,115){\line(0,1){13}}
              \put(363,114){\line(0,1){13}}
              \put(337,60){\line(0,-1){14}}
              \put(361,60){\line(0,-1){14}}
              \put(357,182){\line(0,1){13}}
              \put(337,182){\line(0,1){13}}
              \put(324,60){\framebox(57,55){}}
              \put(323,127){\framebox(57,55){}}
              \put(161,81){\makebox(28,28){A}}
              \put(248,77){\makebox(25,30){B}}
              \put(329,10){\makebox(42,30){AB}}
\end{picture}}

\noindent
{\bf Figure 0 - Tangles and Composition of Tangles} \vspace{3mm}

 {\tt    \setlength{\unitlength}{0.92pt}
\begin{picture}(374,216)
\thinlines    \put(241,44){$cap$}
              \put(329,41){$cup$}
              \put(259,119){$crossings$}
              \put(357,98){\line(0,-1){32}}
              \put(332,66){\line(1,0){24}}
              \put(331,98){\line(0,-1){32}}
              \put(268,96){\line(0,-1){29}}
              \put(245,96){\line(1,0){23}}
              \put(245,67){\line(0,1){28}}
              \put(54,158){\line(1,-1){53}}
              \put(98,160){\line(-2,-3){15}}
              \put(76,127){\line(-3,-5){14}}
              \put(144,158){\line(1,-1){42}}
              \put(182,159){\line(-2,-3){13}}
              \put(163,129){\line(-3,-4){53}}
              \put(107,105){\line(1,-1){19}}
              \put(133,78){\line(1,-1){17}}
              \put(53,159){\line(0,1){32}}
              \put(53,191){\line(1,0){130}}
              \put(182,191){\line(0,-1){32}}
              \put(62,102){\line(0,-1){74}}
              \put(109,57){\line(0,-1){29}}
              \put(150,60){\line(0,-1){34}}
              \put(186,115){\line(0,-1){88}}
              \put(150,26){\line(1,0){37}}
              \put(98,159){\line(1,0){46}}
              \put(62,28){\line(1,0){47}}
\thicklines   \put(12,10){\vector(0,1){196}}
\thinlines    \put(239,183){\line(1,-1){34}}
              \put(237,149){\line(5,4){16}}
              \put(262,168){\line(6,5){14}}
              \put(358,177){\line(-5,-6){29}}
              \put(325,174){\line(6,-5){15}}
              \put(347,154){\line(1,-1){13}}
\end{picture}}

\noindent
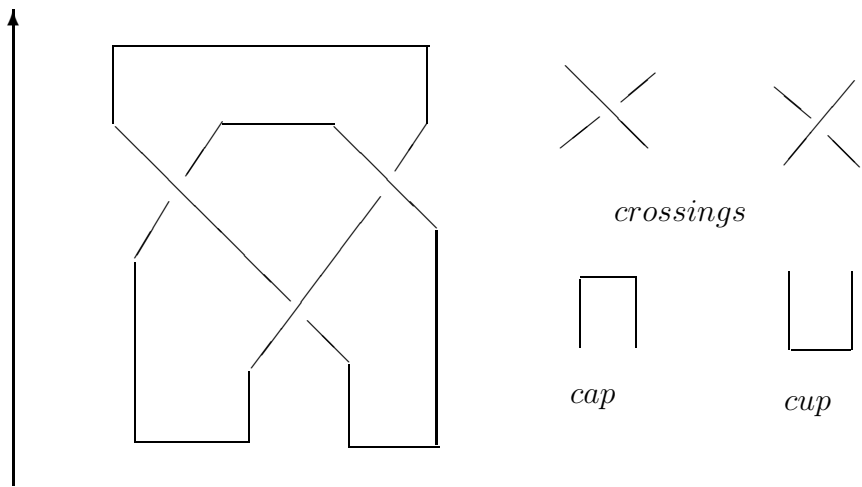
{\bf Figure 1 - Knot as Tangle with Cup and Cap Decomposition}
\vspace{3mm}

Intuitively, a {\em tangle} is a box in three dimensional space
with knotted and linked
string embedded within it and a certain number of strands of that
string emanating from the surface of the box. There are no open
ends of string inside the box. We usually think of some subset of
the strands as {\em inputs} to the tangle and the remaining
strands as the {\em outputs} from the tangle. Usually the inputs
are arranged to be drawn vertically and so that they enter tangle
from below, while the outputs leave the tangle from above. The
tangle itself (within the box) is arranged as nicely as possible
with respect to a vertical direction. This means that a definite
vertical direction is chosen, and that the tangle intersects
planes perpendicular to this direction transversely except for a
finite collection of critical points. These basic critical points
are local maxima and local minima for the space curves inside the
tangle. Two tangles configured with respect to the same box are
{\em ambient isotopic} if there is an  
isotopy in three space carrying one to the other that fixes the input and output
strands of each tangle.  We can {\em compose} two tangles $A$ and
$B$ where the number of output strand of $A$ is equal to the
number of input strands of $B.$ Composition is accomplished by
joining each output strand of $A$ to a corresponding input strand
of $B$. Of course this can be done (in space) in more than one
way. When we use diagrammatic tangles, as we will in this paper,
then the composition operation is naturally well-defined. To have
composition well-defined for spatial tangles we can choose an
ordering of the input and of the output strands of each tangle and
then match them according to this ordering. \vspace{3mm}

Note that the box associated with a tangle or tangle diagram is only a 
delineation of the location of the tangle. The tangle itself consists in the 
woven pattern of strings.  
\vspace{3mm}

We have just given the basic three dimensional description of a
tangle. For the purposes of combinatorial topology and algebra it
is useful to give a modified tangle definition that uses diagrams
instead of embeddings of the tangle into three dimensions. (A
diagram does specify an embedding, but the diagram is not itself
an embedding.) A {\em tangle diagram} is a box in the plane, arranged
parallel to a chosen vertical direction with a left-right ordered sequence
of input strands entering the bottom of the box, and a left-right ordered
sequence of output strands emanating from the top of the box.
Inside the box is a diagram of the tangle represented with
crossings (broken arc indicating the undercrossing line) in the
usual way for knot and links. We assume, as above, that the tangle
is represented so that it is transverse to lines perpendicular to
the vertical except for a finite number of points in the vertical
direction along the tangle. We shall say that the tangle is {\em
well arranged} or {\em Morse} with respect to the vertical direction when these
transversality conditions are met. At the critical points we will
see a local maximum, a local minimum or a crossing in the diagram.
Tangle composition is well-defined (for matching input/output
counts)  since the input and output strands have an ordering (from
left to right for the reader facing the plane on which the tangle
diagram is drawn). Note that the cardinality of the set of input
strands or output strands can be equal to zero. If they are both
zero, then the tangle is simply a knot or link diagram arranged
well with respect to the vertical direction. \vspace{3mm}

The {\em Reidemeister moves} illustrated in 
Figure 2 are a set of moves on diagrams that combinatorially generate isotopy for knots, 
links and and tangles \cite {REID}.  If two tangles are equivalent in three dimensional space,
then corresponding diagrams of these tangles can be obtained one from another, by a sequence
of Reidemeister moves. Each move is confined to the tangle box and keeps the input and 
output strands of the tangle diagram fixed. In illustrating the Reidemeister moves in 
Figure 2 we have shown samples of each type of move. The move 0 is a graphical equivalence 
in the plane that does not change any diagrammatic relations. The move 1 adds or removes a 
twist in the diagram. We have shown one of the two basic examples; the other is obtained by 
switching the crossing in this illustration. Similar remarks apply to obtain other cases of
the moves 2 and 3.
\vspace{3mm}

Two (tangle) diagrams are said to be {\em regularly isotopic} if one can be obtained 
from the other by a sequence of Reidemeister moves of type 0,2,3 (move number 1 is not used 
in regular isotopy).  From now on, all tangles will be tangle diagrams, and we shall 
say that two tangles are {\em equivalent} when they are regularly isotopic.
If we did not insist on arranging our tangle diagrams as Morse diagrams with respect to 
vertical direction, the Reidemeister moves on diagrams would suffice to describe
equivalence of tangles. In order to describe how to move Morse diagrams to one another that
are regularly isotopic, we must add extra moves and rewrite move 0 with respect to the 
vertical. This is illustrated in Figure 3 and Figure 5 and discussed in more detail below.
\vspace{3mm} 

{\tt    \setlength{\unitlength}{0.92pt}
\begin{picture}(357,406)
\thicklines   \put(215,374){\vector(-1,0){25}}
              \put(194,374){\vector(1,0){43}}
              \put(261,374){\line(1,0){86}}
              \put(97,360){\line(4,1){71}}
              \put(125,390){\line(-1,-1){29}}
              \put(52,372){\line(4,1){73}}
              \put(16,293){\framebox(18,19){1}}
              \put(110,308){\vector(1,0){43}}
              \put(134,308){\vector(-1,0){25}}
              \put(149,260){\line(3,2){57}}
              \put(163,341){\line(1,-1){43}}
              \put(72,308){\line(-1,1){21}}
              \put(95,288){\line(-5,4){15}}
              \put(96,316){\line(0,-1){27}}
              \put(38,280){\line(3,2){57}}
              \put(75,203){\line(-3,-5){26}}
              \put(66,178){\line(1,-1){15}}
              \put(47,198){\line(1,-1){12}}
              \put(20,69){\framebox(20,22){3}}
              \put(13,192){\framebox(24,23){2}}
              \put(10,372){\framebox(22,24){0}}
              \put(173,241){\line(0,-1){77}}
              \put(155,241){\line(0,-1){77}}
              \put(57,217){\line(-3,-5){10}}
              \put(71,242){\line(-1,-2){7}}
              \put(46,241){\line(3,-4){28}}
              \put(117,202){\vector(-1,0){25}}
              \put(94,202){\vector(1,0){43}}
              \put(55,130){\line(3,-4){28}}
              \put(80,131){\line(-1,-2){7}}
              \put(66,106){\line(-3,-5){10}}
              \put(83,93){\line(3,-4){28}}
              \put(108,94){\line(-1,-2){7}}
              \put(94,69){\line(-3,-5){10}}
              \put(57,50){\line(3,-4){28}}
              \put(85,53){\line(-1,-2){7}}
              \put(68,28){\line(-3,-5){10}}
              \put(222,131){\line(3,-4){28}}
              \put(247,132){\line(-1,-2){7}}
              \put(233,107){\line(-3,-5){10}}
              \put(199,90){\line(3,-4){28}}
              \put(224,91){\line(-1,-2){7}}
              \put(210,66){\line(-3,-5){10}}
              \put(226,54){\line(3,-4){28}}
              \put(251,53){\line(-1,-2){7}}
              \put(238,31){\line(-3,-5){10}}
              \put(56,89){\line(0,-1){38}}
              \put(108,94){\line(0,1){35}}
              \put(111,55){\line(0,-1){41}}
              \put(200,89){\line(0,1){43}}
              \put(200,48){\line(0,-1){36}}
              \put(250,94){\line(0,-1){43}}
              \put(136,77){\vector(1,0){43}}
              \put(159,77){\vector(-1,0){25}}
\end{picture}}

\noindent
{\bf Figure 2 - Reidemeister Moves}
\vspace{3mm}

If $A$ and $B$ are given tangles, we denote the composition of $A$
and $B$ by $AB$ where the diagram of $A$ is placed below the
diagram of $B$ and the output strands of $A$ are connected to the
input strands of $B$. If the cardinalities of the sets of input
and output strands are zero, then we simple place one tangle below
the other to form the product. \vspace{3mm}

Along with tangle composition, as defined in the previous
paragraph, we also have an operation of {\em product} or {\em
juxtaposition} of tangles. To juxtapose two tangles
$A$ and $B$ simply place their diagrams side by side with $A$ to
the left of $B$ and regard this new diagram as a new tangle whose
inputs are the imputs of $A$ followed by the inputs of $B$, and
whose outputs are the otputs of $A$ followed by the outputs of
$B$. We denote the tangle product of $A$ and $B$ by $A \otimes B$.
\vspace{3mm}

It remains to describe the equivalence relation on tangles that
makes them represent regular isotopy classes of embedded string. For
this purpose it is useful to note that it follows from our
description of the tangle diagrams (See Figure 1) that every tangle is a
composition of {\em elementary tangles} where an elementary tangle
is one of the following list: a {\em cup} (a single minimum --
zero inputs, two outputs), a {\em cap} (a single maximum -- two
inputs, zero outputs), a {\em crossing} (a single local crossing
diagram -- two inputs and two outputs). Figure 5 illustrates
these elementary tangles and the moves on tangles involving
certain compositions of them. These moves include the usual
Reidemeister moves configured with respect to a vertical direction
plus {\em switchback} moves involving the passage of a bit of
string across a maximum or a minimum. We have illustrated these
moves first with the unoriented and then with the different {\em
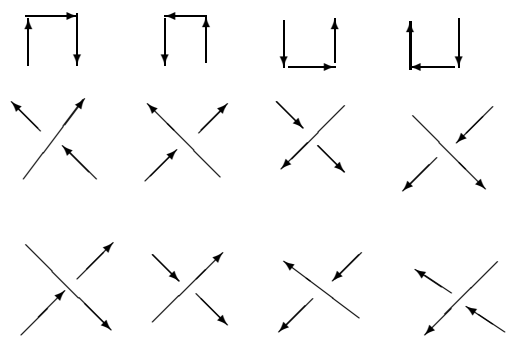
oriented} elementary tangles. Since we consider here {\em regular
isotopy} of tangles, we do not illustrate the first Reidemeister
move which consists in adding or removing a curl (a curl is a
crossing with two of its endpoints connected directly to each
other) from the diagram. Each elementary tangle comes in more than
one flavor due to different possibilities of orientation and
choice of under or over crossing line. Two tangles are said to be
{\em equivalent} if one can be obtained frome the other by a
finite sequence of elementary moves. \vspace{3mm}

In illustrating the elementary moves for oriented tangles, we have
not listed all the cases. In the case of the second Reidemeister
move, we show only the move with reversed orientations on the
lines. In the case of the move of type four, we show only one of
the numerous cases. In most moves there are other cases not illustrated
but obtained from the given illustration by switching one or more crossings.
We leave the full enumeration to the reader.
Note that in this move four a vertical crossing is
exchanged for a horizontal crossing. This relationship allows us
(below) to define the horizontal crossings in terms of the
vertical crossings and the cups and caps. In the type three move
we have illustrated the two main types - all arrows up and two
arrows up, one down. In the discussion below we will show that the
all arrows up or all arrows down move is sufficient to generate
the other type three moves (in the presence of both moves of type
two). \vspace{3mm}

 {\tt    \setlength{\unitlength}{0.92pt}
\begin{picture}(467,310)
\thicklines   \put(80,123){\line(3,-4){12}}
              \put(62,145){\line(4,-5){11}}
              \put(88,148){\line(-3,-5){24}}
              \put(19,32){\framebox(19,21){4}}
              \put(221,136){\framebox(20,22){3}}
              \put(27,138){\framebox(24,23){2}}
              \put(10,246){\framebox(22,24){0}}
              \put(244,42){\line(-4,5){31}}
              \put(233,31){\line(1,1){11}}
              \put(208,12){\line(6,5){14}}
              \put(192,13){\line(0,1){26}}
              \put(192,40){\line(1,0){36}}
              \put(228,38){\line(0,-1){27}}
              \put(72,24){\line(4,-3){13}}
              \put(44,42){\line(5,-3){16}}
              \put(70,80){\line(-3,-4){27}}
              \put(102,42){\line(0,-1){27}}
              \put(66,43){\line(1,0){36}}
              \put(65,16){\line(0,1){26}}
              \put(187,187){\line(0,-1){77}}
              \put(169,187){\line(0,-1){77}}
              \put(71,163){\line(-3,-5){10}}
              \put(85,188){\line(-1,-2){7}}
              \put(60,187){\line(3,-4){28}}
              \put(348,297){\line(0,-1){38}}
              \put(321,259){\line(1,0){25}}
              \put(321,259){\line(0,1){23}}
              \put(298,283){\line(1,0){23}}
              \put(297,283){\line(0,-1){62}}
              \put(192,300){\line(0,-1){89}}
              \put(91,282){\line(0,-1){62}}
              \put(67,284){\line(1,0){23}}
              \put(65,260){\line(0,1){23}}
              \put(40,260){\line(1,0){25}}
              \put(38,298){\line(0,-1){38}}
              \put(141,257){\vector(-1,0){25}}
              \put(118,257){\vector(1,0){43}}
              \put(245,255){\vector(-1,0){25}}
              \put(222,255){\vector(1,0){43}}
              \put(131,148){\vector(-1,0){25}}
              \put(108,148){\vector(1,0){43}}
              \put(146,46){\vector(-1,0){25}}
              \put(123,46){\vector(1,0){43}}
              \put(256,197){\line(3,-4){28}}
              \put(281,198){\line(-1,-2){7}}
              \put(267,173){\line(-3,-5){10}}
              \put(284,160){\line(3,-4){28}}
              \put(309,161){\line(-1,-2){7}}
              \put(295,136){\line(-3,-5){10}}
              \put(258,117){\line(3,-4){28}}
              \put(286,120){\line(-1,-2){7}}
              \put(269,95){\line(-3,-5){10}}
              \put(423,198){\line(3,-4){28}}
              \put(448,199){\line(-1,-2){7}}
              \put(434,174){\line(-3,-5){10}}
              \put(400,157){\line(3,-4){28}}
              \put(425,158){\line(-1,-2){7}}
              \put(411,133){\line(-3,-5){10}}
              \put(427,121){\line(3,-4){28}}
              \put(452,120){\line(-1,-2){7}}
              \put(439,98){\line(-3,-5){10}}
              \put(257,156){\line(0,-1){38}}
              \put(309,161){\line(0,1){35}}
              \put(312,122){\line(0,-1){41}}
              \put(401,156){\line(0,1){43}}
              \put(401,115){\line(0,-1){36}}
              \put(451,161){\line(0,-1){43}}
              \put(337,144){\vector(1,0){43}}
              \put(360,144){\vector(-1,0){25}}
\end{picture}}

\noindent
{\bf Figure 3 - Unoriented Moves on Tangles} \vspace{3mm}

{\tt    \setlength{\unitlength}{0.92pt}
\begin{picture}(224,154)
\thinlines    \put(37,34){\vector(1,1){15}}
              \put(14,11){\vector(1,1){18}}
              \put(16,48){\vector(1,-1){35}}
              \put(208,105){\vector(-1,-1){15}}
              \put(185,84){\vector(-1,-1){14}}
              \put(174,120){\vector(0,1){20}}
              \put(192,121){\vector(-1,0){18}}
              \put(194,141){\vector(0,-1){20}}
              \put(143,123){\vector(0,1){18}}
              \put(124,121){\vector(1,0){19}}
              \put(122,140){\vector(0,-1){19}}
              \put(73,141){\vector(0,-1){19}}
              \put(90,142){\vector(-1,0){17}}
              \put(90,123){\vector(0,1){19}}
              \put(37,143){\vector(0,-1){21}}
              \put(16,142){\vector(1,0){21}}
              \put(17,122){\vector(0,1){19}}
              \put(15,75){\vector(3,4){25}}
              \put(45,75){\vector(-1,1){14}}
              \put(22,95){\vector(-1,1){12}}
              \put(96,76){\vector(-1,1){30}}
              \put(65,74){\vector(1,1){13}}
              \put(87,94){\vector(1,1){12}}
              \put(147,105){\vector(-1,-1){26}}
              \put(119,107){\vector(1,-1){11}}
              \put(136,89){\vector(1,-1){11}}
              \put(175,101){\vector(1,-1){30}}
              \put(68,16){\vector(1,1){29}}
              \put(68,44){\vector(1,-1){11}}
              \put(86,28){\vector(1,-1){13}}
              \put(153,18){\vector(-4,3){31}}
              \put(154,45){\vector(-1,-1){12}}
              \put(134,26){\vector(-1,-1){14}}
              \put(210,41){\vector(-1,-1){30}}
              \put(213,12){\vector(-3,2){16}}
              \put(191,28){\vector(-3,2){15}}
\end{picture}}

\noindent
{\bf Figure 4 - Oriented Elementary Tangles} \vspace{3mm}

 {\tt    \setlength{\unitlength}{0.92pt}
\begin{picture}(404,339)
\thinlines    \put(352,148){\vector(2,3){9}}
              \put(169,151){\vector(2,3){9}}
              \put(359,121){\vector(-2,-3){11}}
              \put(381,148){\vector(-3,-4){11}}
              \put(347,181){\vector(1,-1){33}}
              \put(240,139){\vector(1,-1){36}}
              \put(250,155){\vector(-2,-3){11}}
              \put(266,180){\vector(-2,-3){8}}
              \put(293,106){\vector(-2,3){52}}
              \put(273,152){\vector(1,2){15}}
              \put(326,103){\vector(1,2){19}}
              \put(375,104){\vector(-2,3){52}}
              \put(245,102){\vector(1,2){8}}
              \put(259,126){\vector(1,2){8}}
              \put(369,172){\vector(1,2){6}}
              \put(288,129){\makebox(40,28){=}}
              \put(106,133){\makebox(40,28){=}}
              \put(194,149){\vector(-3,4){28}}
              \put(165,105){\vector(2,3){12}}
              \put(185,172){\vector(1,2){6}}
              \put(185,134){\vector(2,3){10}}
              \put(78,172){\vector(2,3){9}}
              \put(61,146){\vector(2,3){11}}
              \put(89,109){\vector(-3,4){27}}
              \put(77,130){\vector(1,2){8}}
              \put(63,106){\vector(1,2){8}}
              \put(193,108){\vector(-2,3){52}}
              \put(144,107){\vector(1,2){19}}
              \put(91,156){\vector(1,2){15}}
              \put(113,109){\vector(-2,3){52}}
              \put(19,192){\vector(0,1){73}}
              \put(20,266){\vector(1,0){27}}
              \put(47,266){\vector(0,-1){29}}
              \put(48,237){\vector(1,0){22}}
              \put(70,236){\vector(0,1){93}}
              \put(104,199){\vector(0,1){128}}
              \put(76,242){\makebox(18,22){=}}
              \put(187,272){\vector(-1,0){22}}
              \put(166,240){\vector(-1,0){17}}
              \put(148,241){\vector(0,1){78}}
              \put(189,201){\vector(0,1){71}}
              \put(166,269){\vector(0,-1){29}}
              \put(200,249){\makebox(18,22){=}}
              \put(228,201){\vector(0,1){128}}
              \put(315,251){\makebox(12,14){=}}
              \put(304,209){\vector(0,1){101}}
              \put(274,311){\vector(0,-1){103}}
              \put(367,232){\vector(-1,1){23}}
              \put(392,207){\vector(-1,1){19}}
              \put(392,258){\vector(-3,-4){39}}
              \put(375,292){\vector(3,4){17}}
              \put(345,257){\vector(3,4){18}}
              \put(352,312){\vector(3,-4){40}}
              \put(106,10){\vector(0,1){38}}
              \put(106,48){\vector(1,0){33}}
              \put(139,47){\vector(0,-1){36}}
              \put(122,12){\vector(-1,1){13}}
              \put(101,30){\vector(-1,1){12}}
              \put(88,43){\vector(1,1){34}}
              \put(223,48){\vector(0,-1){36}}
              \put(190,49){\vector(1,0){33}}
              \put(190,11){\vector(0,1){38}}
              \put(227,32){\vector(1,1){15}}
              \put(242,47){\vector(-4,3){35}}
              \put(204,14){\vector(1,1){14}}
              \put(145,17){\makebox(40,28){=}}
              \put(10,284){\framebox(19,20){0}}
              \put(246,264){\framebox(17,21){2}}
              \put(20,136){\framebox(22,23){3}}
              \put(37,29){\framebox(21,25){4}}
\end{picture}}

\noindent
{\bf Figure 5 - Representative Moves on Oriented Tangles} \vspace{3mm}

In considering the oriented moves on tangles we see that there are
two basic types of Reidemeister $2$ move. The first ( $2_{A}$) has
both strings oriented in parallel to each other. The second
($2_{B}$) has the strings oriented in opposite directions.
Similarly there are two basic types of third Reidimeister move
that we denote by $3_{C}$ and $3_{NC}$ where the former has a
cyclic triangle and the latter does not. See Figure 5. It turns
out that the cyclic type three move is a consequence of the
reverse oriented two move and the non-cyclic type three move. The
proof of this statement is illustrated in Figure 6. \vspace{3mm}

{\tt    \setlength{\unitlength}{0.92pt}
\begin{picture}(405,479)
\thinlines    \put(88,80){\vector(1,0){13}}
              \put(15,79){\vector(1,0){64}}
              \put(21,164){\vector(3,-4){101}}
              \put(331,77){\vector(-1,-2){29}}
              \put(355,120){\vector(-1,-2){16}}
              \put(374,154){\vector(-1,-2){13}}
              \put(99,161){\vector(-2,-3){20}}
              \put(31,56){\vector(-2,-3){19}}
              \put(43,75){\vector(-2,-3){6}}
              \put(60,102){\vector(-2,-3){11}}
              \put(73,124){\vector(-2,-3){7}}
              \put(305,340){\vector(-1,-2){14}}
              \put(316,361){\vector(-1,-2){6}}
              \put(325,381){\vector(-2,-3){7}}
              \put(332,398){\vector(-1,-2){5}}
              \put(366,460){\vector(-1,-2){25}}
              \put(63,388){\vector(-1,-2){14}}
              \put(38,346){\vector(-1,-2){20}}
              \put(102,456){\vector(-1,-2){27}}
              \put(45,141){\vector(2,3){17}}
              \put(31,125){\vector(3,4){9}}
              \put(44,124){\vector(-1,0){16}}
              \put(101,126){\vector(-1,0){44}}
              \put(384,367){\vector(0,1){17}}
              \put(342,384){\vector(-1,0){52}}
              \put(383,383){\vector(-1,0){27}}
              \put(371,367){\vector(1,0){13}}
              \put(287,366){\vector(1,0){69}}
              \put(292,459){\vector(3,-4){101}}
              \put(102,81){\vector(0,1){45}}
              \put(68,26){\vector(-1,1){53}}
              \put(339,314){\vector(-1,1){53}}
              \put(290,382){\vector(1,2){19}}
              \put(320,436){\vector(1,2){16}}
              \put(379,103){\vector(-1,1){52}}
              \put(333,15){\vector(1,2){17}}
              \put(362,69){\vector(1,2){17}}
              \put(286,146){\vector(3,-4){101}}
              \put(20,457){\vector(3,-4){101}}
              \put(51,436){\vector(1,2){16}}
              \put(21,382){\vector(1,2){19}}
              \put(76,308){\vector(-3,4){55}}
\thicklines   \put(124,387){\vector(1,0){123}}
              \put(267,294){\vector(-4,-3){137}}
              \put(141,93){\vector(1,0){127}}
\thinlines    \put(130,391){\framebox(113,26){Reverse 2}}
              \put(211,199){\framebox(126,41){Non Cyclic 3}}
              \put(148,98){\framebox(117,32){Reverse 2}}
\end{picture}}

\noindent
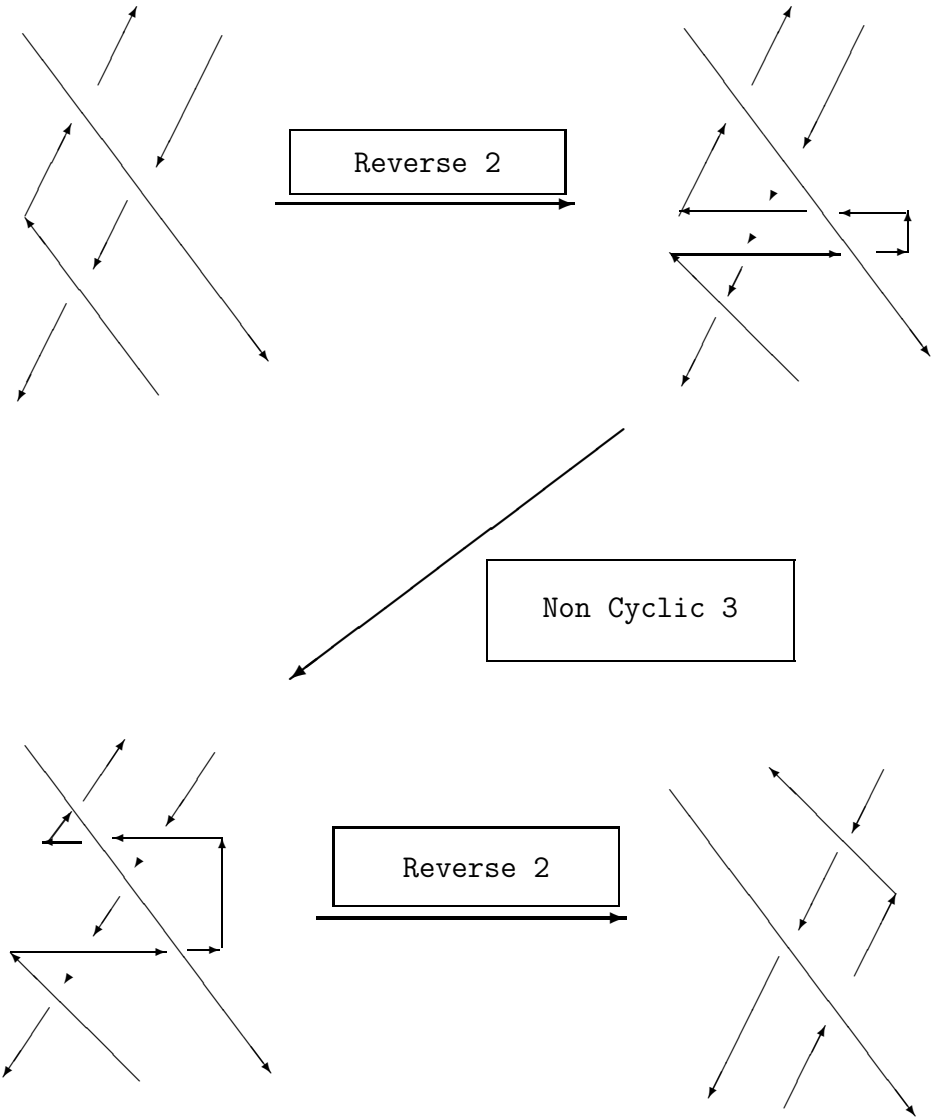
{\bf Figure 6 - Reverse Type Two Move Plus Non Cyclic Type three
Implies Cyclic Type Three Move} \vspace{3mm}

This reduction of the number of type three moves makes working
with the tangle category easier. It is still the case that one
must verify invariance under both varieties of type two move for
any functor on $Tang$. Note also the switchback moves shown in
Figure 5.  These moves give necessary relations between the
crossings and the maxima and minima. \vspace{3mm}

In order to make $Tang$ a category, we define the {\em objects} of
$Tang$ to be ordered lists of signs ($+1$ or $-1$), including the
empty list. Thus the objects of $Tang$ are the lists $$<>,
<\epsilon_{1}>, <\epsilon_{1},\epsilon_{2}>, ...,
<\epsilon_{1},\epsilon_{2},...,\epsilon_{n}>, ...$$ We let
$\epsilon$ denote the list
$<\epsilon_{1},\epsilon_{2},...,\epsilon_{n}>$ with $0$ standing
for the empty list. We let $||\epsilon||$ denote the number of
entries in the list $\epsilon$. A tangle is a morphism from
$\epsilon$ to $\epsilon'$ where $||\epsilon||$ and $||\epsilon'||$
are the number of input and output strands of this tangle,
respectively. Each sign corresponds to a single input or output
strand of the tangle. If the sign is positive, then the strand is
oriented upwards. If the sign is negative, then the strand is
oriented downwards. Two morphisms are equal if they are equivalent
as tangles. It is clear that this assignment makes $Tang$ into a
category, where composition of tangles corresponds to composition
of morphisms in the category. The juxtapositon or tensor operation
on tangles makes $Tang$ into an associative tensor category with
$\epsilon \otimes \epsilon'$ the list obtained by appending the
second list of signs to the first list of signs. Note that the
empty list is an identity element for the tensor product.
\vspace{3mm}

\subsection{Flat Tangles} We now discuss the simpler category
$Flat$ of flat tangles.  A {\em flat tangle} is exactly the same
as an ordinary tangle diagram except that the crossings in the
flat tangle have no over or under specification associated with
them. Thus the generating morphisms of $Flat$ are the oriented
cups and caps and the various orientations of the flat crossing.
All the generalized Reidemeister moves hold in $Flat$ with no
stipulations about under and over crossings. This means that there
is much more freedom to perfom moves in $Flat$. See Figure 7 for
sample illustrations of some of the moves in $Flat.$ Note that the
crossing in $Flat$ has all the formal properties of a permutation.
In fact, the category $Flat$ naturally contains copies of all the
symmetric groups with the symmetric group on $n$ letters occuring
as the flat tangles on $n$ upward oriented strands that do not use
any cups or caps. \vspace{3mm}

{\tt    \setlength{\unitlength}{0.92pt} \begin{picture}(408,341)
\thinlines    \put(202,11){\vector(1,1){38}}
\put(118,14){\vector(-1,1){30}} \put(377,147){\vector(-3,-4){32}}
\put(324,109){\vector(2,3){50}} \put(267,186){\vector(-2,-3){27}}
\put(244,110){\vector(2,3){45}} \put(168,107){\vector(2,3){26}}
\put(145,109){\vector(2,3){49}} \put(62,147){\vector(3,4){28}}
\put(64,112){\vector(2,3){46}} \put(397,214){\vector(-1,1){46}}
\put(350,261){\vector(3,4){38}} \put(344,181){\vector(1,-1){33}}
\put(239,144){\vector(1,-1){36}} \put(293,108){\vector(-2,3){52}}
\put(375,106){\vector(-2,3){52}} \put(288,131){\makebox(40,28){=}}
\put(106,135){\makebox(40,28){=}} \put(192,147){\vector(-3,4){28}}
\put(89,111){\vector(-3,4){27}} \put(193,110){\vector(-2,3){52}}
\put(113,111){\vector(-2,3){52}} \put(19,194){\vector(0,1){73}}
\put(20,268){\vector(1,0){27}} \put(47,268){\vector(0,-1){29}}
\put(48,239){\vector(1,0){22}} \put(70,238){\vector(0,1){93}}
\put(104,201){\vector(0,1){128}} \put(76,244){\makebox(18,22){=}}
\put(187,274){\vector(-1,0){22}} \put(166,242){\vector(-1,0){17}}
\put(148,243){\vector(0,1){78}} \put(189,203){\vector(0,1){71}}
\put(166,271){\vector(0,-1){29}} \put(200,251){\makebox(18,22){=}}
\put(228,203){\vector(0,1){128}} \put(315,253){\makebox(12,14){=}}
\put(304,211){\vector(0,1){101}} \put(274,313){\vector(0,-1){103}}
\put(392,260){\vector(-3,-4){39}} \put(352,314){\vector(3,-4){40}}
\put(106,12){\vector(0,1){38}} \put(106,50){\vector(1,0){33}}
\put(139,49){\vector(0,-1){36}} \put(88,45){\vector(1,1){34}}
\put(223,50){\vector(0,-1){36}} \put(190,51){\vector(1,0){33}}
\put(190,13){\vector(0,1){38}} \put(242,49){\vector(-4,3){35}}
\put(145,19){\makebox(40,28){=}} \put(10,286){\framebox(19,20){0}}
\put(246,266){\framebox(17,21){2}}
\put(20,138){\framebox(22,23){3}} \put(37,31){\framebox(21,25){4}}
\end{picture}}

\noindent
{\bf Figure 7 - Representative Oriented Moves in $Flat$} \vspace{3mm}

\section{The Category of an Oriented Quantum Algebra} Let $A$ be
an oriented quantum algebra. In this section we build a category
$Cat(A)$ associated to $A$ by decorating morphisms of the category $Flat$ of
flat tangles with "beads" from the quantum algebra. The category
$Cat(A)$ has built in bead sliding rules that allow the reduction of
individual strings to pure algebra.  We will construct a functor
from the tangle category to the category of the quantum algebra
and  show this functor is well-defined. In the course of this
construction the reasons for our definition of oriented quantum
algebra will become transparent. We will discuss the structure of
evaluations on $Cat(A)$ that give rise to invariants of links and
tangles. \vspace{3mm}

First note that given any algebra $A$ there is a category $C(A)$,
associated with this algebra with one object $[+]$ and a morphism
$a$ for every element $a$ in the algebra (we use the same symbol
for the element and its corresponding morphism). We simply declare
that $a:[+] \longrightarrow [+]$ and that composition of morphisms
corresponds to the multiplication of elements in the algebra.
Since $A$ is an algebra, we can also add morphisms of $C(A)$ in
correspondence to the addition of elements of $A$. \vspace{3mm}

The category $Cat(A)$ will be constructed by thinking of the morphisms in 
$C(A)$ as labelled arrows and generalizing them to labelled flat oriented 
diagrams from the flat tangle category. In this way, we see at once that 
we must deal with down arrows as well as up arrows, and since down arrows in 
$Flat$ terminate in the object $[-]$, we must enlarge the list of objects to 
those generated by the empty object $[~]$ and the objects $[+]$ and $[-]$ from 
$Flat$.
\vspace{3mm}

Using the same object structure that we described for the tangle
categories, we can regard the basic element for
$C(A^{\otimes^{n}})$ as the object 

$$[n] = [+] \otimes [+] \otimes... \otimes [+]$$ 

\noindent ($n$ factors). Morphisms in
$C_{n}(A)=C(A^{\otimes^{n}})$ are sums of morphisms that correspond
to $n$-fold tensor products of elements of $A$. However,by the motivation above,
we also need morphisms whose objects are negative signed lists. The
simplest such object is $[-]$, and a morphism from $[-1]$ to
$[-1]$ is a downward pointing arrow in $Flat$. We can decorate
this arrow with an element $a$ of the algebra $A$. As a morphism
we denote it by $a^{t}:[-1] \longrightarrow [-1]$.  Formally, the
morphisms $a:[+] \longrightarrow [+]$ and $a^{t}$ determine each
other. Informally, $a^{t}$ is just what you see if you reverse the
sense of external vertical direction for the morphisms of
$Cat(A)$. \vspace{3mm}

We can bundle all of these categories together into one tensor
category $C_{\infty}(A)$ with morphisms in $1-1$ correspondence
with the elements in arbitrary-fold tensor products of $A$ and
objects in $1-1$ correspondence with finite sequences of signs.
\vspace{3mm}

 {\tt    \setlength{\unitlength}{0.92pt}
\begin{picture}(371,487)
\thinlines    \put(312,231){$a^{t}$}
              \put(293,463){$[-]$}
              \put(296,326){$[-]$}
              \put(108,329){$[+]$}
              \put(109,469){$[+]$}
              \put(108,60){\makebox(27,28){ab}}
              \put(295,55){\makebox(27,28){ab}}
              \put(210,33){\makebox(28,26){b}}
              \put(207,58){\makebox(28,25){a}}
              \put(10,60){\makebox(24,34){b}}
              \put(10,40){\makebox(24,25){a}}
              \put(30,159){\framebox(77,150){}}
              \put(234,165){\framebox(74,140){}}
              \put(244,223){\makebox(18,21){a}}
              \put(45,220){\makebox(17,23){a}}
              \put(286,70){\circle*{12}}
              \put(285,131){\vector(0,-1){120}}
              \put(247,69){\circle*{12}}
              \put(246,130){\vector(0,-1){120}}
              \put(96,11){\vector(0,1){125}}
              \put(97,75){\circle*{12}}
              \put(40,13){\vector(0,1){125}}
              \put(41,77){\circle*{12}}
              \put(274,295){\vector(0,-1){120}}
              \put(45,404){\circle*{12}}
              \put(247,47){\circle*{12}}
              \put(41,53){\circle*{12}}
              \put(275,234){\circle*{12}}
              \put(74,172){\vector(0,1){125}}
              \put(44,341){\vector(0,1){125}}
              \put(74,232){\circle*{12}}
              \put(20,392){\makebox(17,23){a}}
              \put(97,392){\makebox(17,23){a}}
              \put(119,340){\vector(0,1){125}}
              \put(121,404){\circle*{12}}
              \put(64,400){$=$}
              \put(234,456){\vector(0,-1){120}}
              \put(303,457){\vector(0,-1){120}}
              \put(304,398){\circle*{12}}
              \put(235,396){\circle*{12}}
              \put(208,384){\makebox(17,23){a}}
              \put(277,385){\makebox(17,23){a}}
              \put(251,394){$=$}
              \put(113,228){$a$}
\end{picture}}

\noindent
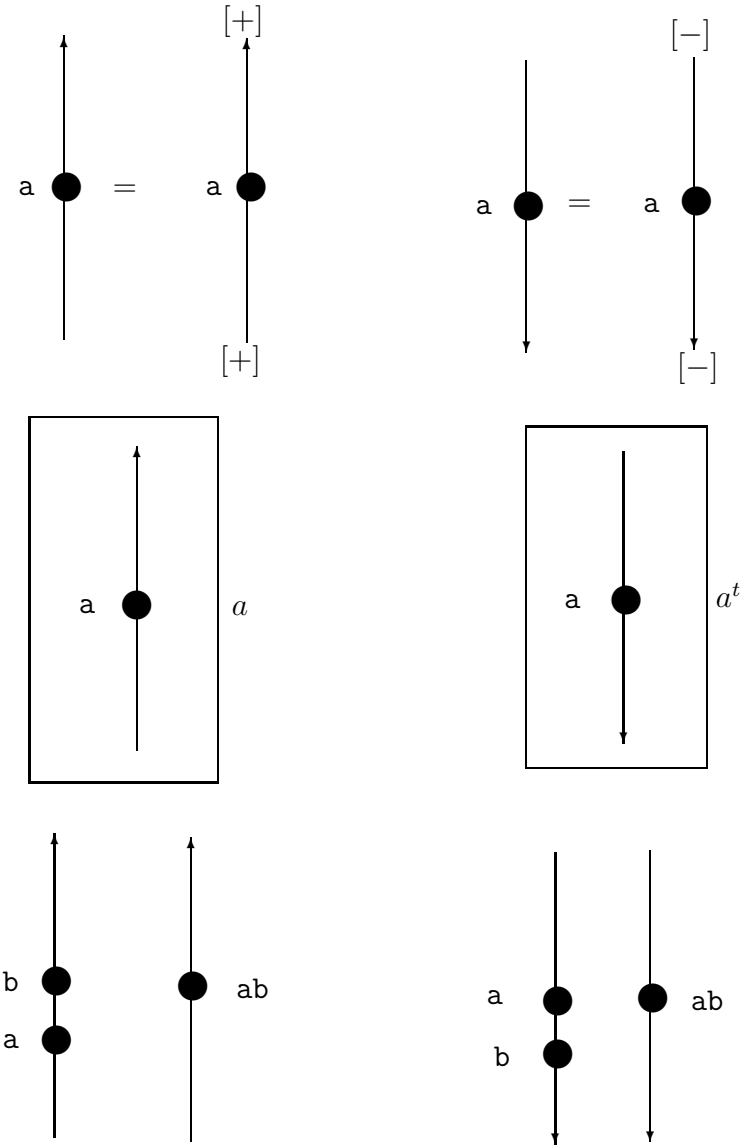
{\bf Figure 8 - Decorated Single Arrows and the Transpose of a Morphism}
\vspace{3mm}

For a quantum algebra $A$, we now extend this tensor category
$C_{\infty}(A)$ to encompass both the automorphisms $U$ and $D$ of the quantum
algebra and the structure of the flat tangle category $Flat.$ We
accomplish this by interpreting morphisms in $C_{1}(A)$ as single
vertical lines (that is, as flat $1-1$ tangles) decorated by a
label that corresponds to the algebra element that gives this
morphism. In our notation, this decoration is a round node (a
"bead") on the line that is usually accompanied by a text label in
the plane next to the bead and disjoint from the line. Since the
algebra element may be a sum of elements, we may use a sum of
such labelled  vertical segments. An unlabelled single segment
denotes the identity morphism in $C_{1}(A).$ Similarly, a morphism
in $C_{n}(A)$ is denoted by a labelling of a set of $n$ parallel
vertical segments. If the segments (from left to right) are
labelled $a_{1}$, $a_{2}$, ... $a_{n}$, then this labelled bundle
of segments is the morphism in $C_{n}(A)$ that corresponds to the
tensor product

$$a_{1} \otimes a_{2} \otimes  ... \otimes a_{n}.$$

\noindent So far, we have indicated how to represent the category
$C_{\infty}(A)$ via labelling of the identity $n$-tangles in
$Flat.$ We now extend this to a category we call {\em the
category of the oriented quantum algebra $A$} and denote by
$Cat(A).$ The objects in $Cat(A)$ are identical to the objects in
$Flat$. The morphisms in $Cat(A)$ are flat tangles that have been
decorated on some of their vertical segments by elements of the
algebra $A$. (Thus all the elements of $C_{\infty}(A)$ are
reperesented in $Cat(A)$ by decorations of the identity tangles.)
The relations on generating morphisms for $Flat$ still hold in
$Cat(A)$. In addition we have extra relations for the cups and
caps as illustrated in Figure 10.  The basic form of this relation
is:

$$ (a \otimes 1_{V}) \circ Cup =  (1_{V} \otimes \tau(a)) \circ Cup,
$$

$$ Cap \circ (1_{V} \otimes a) =  Cap \circ ( \tau'(a) \otimes
1_{V})$$

\noindent where $\tau, \tau': A \longrightarrow A$ are automorphisms of the
algebra $A.$

\noindent Here we have not specified the orientations on the $Cup$
or the $Cap.$ Up to the choice of the automorphism $\tau$ or $\tau'$, these
relations are independent of the choice of orientation, and apply
to any element $a$ of the algebra $A$. The topological relations
for $Cup$ and $Cap$ that hold in $Flat$ are extended to $Cat(A)$. 
Once we orient the cups and caps we can specify the choice of
automorphism $\tau$ or $\tau'$ from two possibilities that we call $U$ ("up")
and $D$ ("down").  Caps with clockwise orientation and cups with
counterclockwise orientation receive $U$. Caps with
counterclockwise orientation and cups with clockwise orientation
reveive $D.$  See Figures 9 and 10. \vspace{3mm}

The upshot of these relations of the $Cup$ and $Cap$ with the
automorphisms $U$ and $D$ of $A$ is that $U$ and $D$ have
diagrammatic interpretations as shown in Figure 10. One way to
think about this diagrammatic interpretation is that $U(a)$ is
obtained from the upward pointing diagram for the morphism $a$ by
turning this diagram upside down and running a vertical line
upward that turns through a cap (on the left) and connects to the
top part of the inverted $a$. Then a $Cup$ is connected to the
bottom part and continues upward to form a globally upward
pointing morphism that represents $U(a)$. For $D(a)$ the same
diagram is used but all the arrows are reversed. \vspace{3mm}

 {\tt    \setlength{\unitlength}{0.92pt}
\begin{picture}(405,343)
\thinlines    \put(329,94){$D^{-1}a$}
              \put(332,251){$U^{-1}a$}
              \put(288,94){\makebox(18,22){=}}
              \put(127,98){\makebox(18,22){=}}
              \put(71,109){\makebox(16,17){a}}
              \put(231,109){\makebox(16,17){a}}
              \put(319,180){\vector(0,-1){167}}
              \put(157,177){\vector(0,-1){167}}
              \put(319,97){\circle*{8}}
              \put(252,116){\circle*{8}}
              \put(93,118){\circle*{8}}
              \put(157,100){\circle*{8}}
              \put(275,134){\vector(0,-1){118}}
              \put(253,134){\vector(1,0){23}}
              \put(251,103){\vector(0,1){31}}
              \put(227,103){\vector(1,0){24}}
              \put(228,176){\vector(0,-1){73}}
              \put(63,132){\vector(0,-1){117}}
              \put(93,136){\vector(-1,0){30}}
              \put(93,104){\vector(0,1){33}}
              \put(117,104){\vector(-1,0){25}}
              \put(118,179){\vector(0,-1){75}}
              \put(159,237){\makebox(44,34){$Ua$}}
              \put(261,250){\makebox(16,17){a}}
              \put(323,255){\circle*{8}}
              \put(322,202){\vector(0,1){128}}
              \put(294,250){\makebox(18,22){=}}
              \put(260,271){\vector(0,-1){29}}
              \put(261,258){\circle*{8}}
              \put(283,202){\vector(0,1){71}}
              \put(242,242){\vector(0,1){78}}
              \put(260,241){\vector(-1,0){17}}
              \put(281,273){\vector(-1,0){22}}
              \put(124,246){\makebox(18,22){=}}
              \put(152,203){\vector(0,1){128}}
              \put(153,256){\circle*{8}}
              \put(77,248){\makebox(16,17){a}}
              \put(19,251){\makebox(16,17){a}}
              \put(14,261){\circle*{8}}
              \put(13,194){\vector(0,1){128}}
              \put(96,256){\circle*{8}}
              \put(118,240){\vector(0,1){93}}
              \put(96,241){\vector(1,0){22}}
              \put(95,270){\vector(0,-1){29}}
              \put(68,270){\vector(1,0){27}}
              \put(67,196){\vector(0,1){73}}
              \put(165,85){\makebox(35,29){$Da$}}
\end{picture}}

\noindent
{\bf Figure 9 - Diagrams for the Automorphisms $U$ and $D$} \vspace{3mm}

{\tt    \setlength{\unitlength}{0.92pt} \begin{picture}(469,507)
\thinlines    \put(258,10){\framebox(200,86){}}
\put(188,228){\framebox(0,0){}} \put(258,113){\framebox(200,86){}}
\put(345,47){\makebox(18,22){=}} \put(344,147){\makebox(18,22){=}}
\put(274,44){\makebox(18,24){a}}
\put(421,43){\makebox(35,29){$Da$}}
\put(361,147){\makebox(35,29){$Da$}}
\put(325,148){\makebox(18,24){a}} \put(412,74){\vector(0,-1){34}}
\put(412,40){\vector(-1,0){36}} \put(376,40){\vector(0,1){33}}
\put(438,144){\vector(0,1){34}} \put(437,179){\vector(-1,0){36}}
\put(401,177){\vector(0,-1){35}} \put(298,57){\circle*{8}}
\put(413,58){\circle*{8}} \put(402,161){\circle*{8}}
\put(318,161){\circle*{8}} \put(298,41){\vector(0,1){33}}
\put(334,41){\vector(-1,0){36}} \put(334,75){\vector(0,-1){34}}
\put(281,178){\vector(0,-1){35}} \put(317,180){\vector(-1,0){36}}
\put(318,145){\vector(0,1){34}} \thicklines  
\put(197,105){\vector(1,0){41}} \thinlines   
\put(262,230){\framebox(196,86){}}
\put(262,358){\framebox(197,102){}}
\put(325,253){\makebox(18,22){=}}
\put(407,247){\makebox(44,34){$Ua$}}
\put(286,257){\makebox(16,17){a}} \put(400,265){\circle*{8}}
\put(279,267){\circle*{8}} \put(362,283){\vector(0,-1){36}}
\put(361,247){\vector(1,0){38}} \put(399,249){\vector(0,1){34}}
\put(316,250){\vector(0,1){34}} \put(278,248){\vector(1,0){38}}
\put(279,284){\vector(0,-1){36}} \put(12,205){\vector(0,1){73}}
\put(13,279){\vector(1,0){27}} \put(40,279){\vector(0,-1){29}}
\put(41,250){\vector(1,0){22}} \put(63,249){\vector(0,1){93}}
\put(41,265){\circle*{8}} \put(22,257){\makebox(16,17){a}}
\put(99,266){\circle*{8}} \put(98,213){\vector(0,1){128}}
\put(69,285){\makebox(18,22){=}}
\put(105,247){\makebox(44,34){$Ua$}} \put(206,273){\circle*{8}}
\put(206,254){\vector(0,1){93}} \put(184,254){\vector(1,0){22}}
\put(183,283){\vector(0,-1){29}} \put(156,283){\vector(1,0){27}}
\put(156,210){\vector(0,1){73}} \put(122,286){\makebox(18,22){=}}
\thicklines   \put(217,271){\vector(1,0){35}} \thinlines   
\put(337,400){\makebox(18,22){=}}
\put(315,404){\makebox(16,17){a}}
\put(356,395){\makebox(44,34){$Ua$}} \put(404,414){\circle*{8}}
\put(310,412){\circle*{8}} \put(404,395){\vector(0,1){36}}
\put(405,431){\vector(1,0){34}} \put(439,430){\vector(0,-1){38}}
\put(309,428){\vector(0,-1){38}} \put(275,429){\vector(1,0){34}}
\put(274,393){\vector(0,1){36}} \thicklines  
\put(224,421){\vector(1,0){35}} \thinlines   
\put(125,413){\makebox(18,22){=}} \put(163,360){\vector(0,1){73}}
\put(163,433){\vector(1,0){27}} \put(190,433){\vector(0,-1){29}}
\put(191,404){\vector(1,0){22}} \put(213,404){\vector(0,1){93}}
\put(163,413){\circle*{8}} \put(102,114){\makebox(18,22){=}}
\put(44,125){\makebox(16,17){a}} \put(130,193){\vector(0,-1){167}}
\put(66,135){\circle*{8}} \put(130,116){\circle*{8}}
\put(35,152){\vector(0,-1){117}} \put(66,152){\vector(-1,0){30}}
\put(66,120){\vector(0,1){33}} \put(90,120){\vector(-1,0){25}}
\put(91,195){\vector(0,-1){75}}
\put(113,386){\makebox(44,34){$Ua$}}
\put(75,413){\makebox(18,22){=}} \put(105,363){\vector(0,1){128}}
\put(106,416){\circle*{8}} \put(29,407){\makebox(16,17){a}}
\put(48,415){\circle*{8}} \put(70,399){\vector(0,1){93}}
\put(48,400){\vector(1,0){22}} \put(47,429){\vector(0,-1){29}}
\put(20,429){\vector(1,0){27}} \put(19,355){\vector(0,1){73}}
\put(144,97){\makebox(35,29){$Da$}} \end{picture}}

\noindent
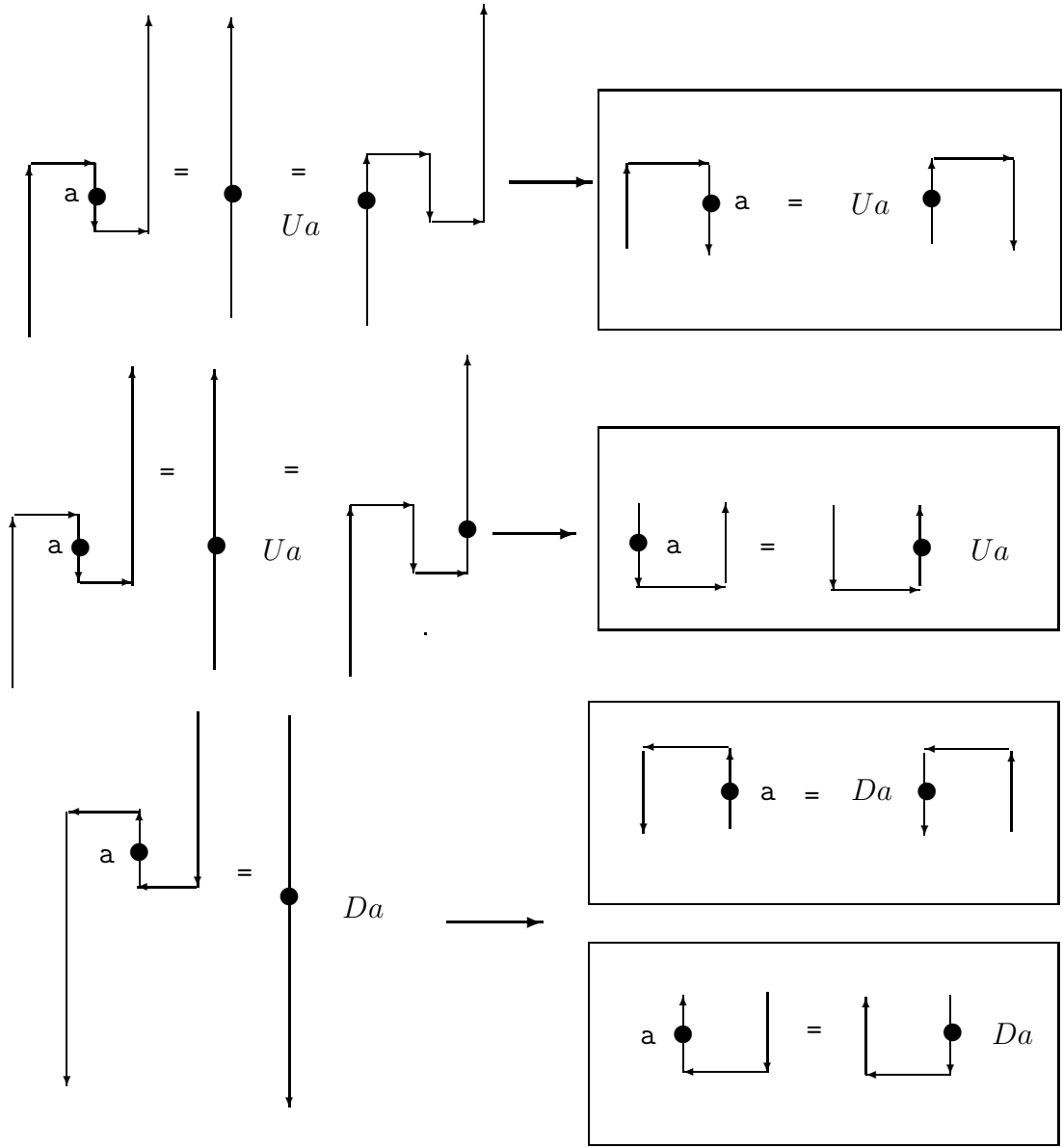
{\bf Figure 10 - Bead Sliding and the Automorphisms $U$ and $D$}
\vspace{3mm}

\noindent The reader may wonder if it is necessary to have two
distinct automorphisms $U$ and $D.$ In fact, we shall see that
there are cases where it is most convenient to have $U=D$ and
other cases where it is most natural to take $D$ to be the
identity mapping while $U$ is non-trivial. The latter occurs
in representing a quasi-triangular ribbon Hopf algebra, as we will
see later in section 6. \vspace{3mm}

It is interesting to note that the morphisms $UD$ and $DU$ are
both of the form $UD(a) = DU(a) = GaG^{-1}$ where $G$ is the flat
curl morphism illustrated in Figure 11. In the case of a ribbon
Hopf algebra it is convenient to interpret $G$ as a certain
grouplike element in the algebra itself. \vspace{3mm}

{\tt    \setlength{\unitlength}{0.92pt} \begin{picture}(341,349)
\thinlines    \put(31,10){\makebox(300,56){$UDa = DUa =
GaG^{-1}$}} \put(80,160){\vector(0,1){18}}
\put(35,142){\vector(0,1){18}} \put(138,92){\vector(0,-1){14}}
\put(94,122){\vector(0,-1){14}} \put(250,93){\vector(0,-1){14}}
\put(57,142){\vector(0,-1){14}} \put(324,182){\vector(0,-1){14}}
\put(131,231){\vector(0,1){18}} \put(242,231){\vector(0,1){18}}
\put(316,319){\vector(0,1){18}} \put(27,253){\vector(0,-1){14}}
\put(213,181){\vector(0,-1){14}} \put(20,94){\vector(0,-1){14}}
\put(72,270){\vector(0,-1){14}} \put(12,291){\vector(0,1){18}}
\put(49,273){\vector(0,1){18}} \put(205,318){\vector(0,1){18}}
\put(86,320){\vector(0,1){18}} \put(36,135){\makebox(17,19){a}}
\put(222,138){\makebox(22,19){=}}
\put(111,138){\makebox(22,19){=}} \put(310,127){\line(0,1){12}}
\put(265,147){\line(0,1){15}} \put(310,140){\line(-3,-1){59}}
\put(324,167){\line(-3,-1){59}} \put(310,126){\line(-1,0){22}}
\put(287,162){\line(-1,0){22}} \put(287,161){\line(0,-1){35}}
\put(288,143){\circle*{10}} \put(324,167){\line(0,1){40}}
\put(250,80){\line(0,1){40}} \put(138,78){\line(0,1){40}}
\put(213,167){\line(0,1){40}} \put(139,119){\line(1,1){59}}
\put(154,108){\line(1,1){59}} \put(177,143){\circle*{10}}
\put(176,161){\line(0,-1){35}} \put(176,162){\line(-1,0){22}}
\put(199,126){\line(-1,0){22}} \put(154,162){\line(0,-1){53}}
\put(199,179){\line(0,-1){53}} \put(94,209){\line(0,-1){99}}
\put(20,180){\line(0,-1){99}} \put(21,180){\line(1,0){59}}
\put(35,109){\line(1,0){59}} \put(80,180){\line(0,-1){53}}
\put(35,163){\line(0,-1){53}} \put(80,127){\line(-1,0){22}}
\put(57,163){\line(-1,0){22}} \put(57,163){\line(0,-1){35}}
\put(59,146){\circle*{10}} \put(49,272){\circle*{10}}
\put(49,292){\line(0,-1){35}} \put(49,293){\line(-1,0){22}}
\put(72,257){\line(-1,0){22}} \put(27,293){\line(0,-1){53}}
\put(72,310){\line(0,-1){53}} \put(27,239){\line(1,0){59}}
\put(13,310){\line(1,0){59}} \put(12,310){\line(0,-1){99}}
\put(86,339){\line(0,-1){99}} \put(191,309){\line(0,-1){53}}
\put(146,292){\line(0,-1){53}} \put(191,256){\line(-1,0){22}}
\put(168,292){\line(-1,0){22}} \put(168,291){\line(0,-1){35}}
\put(169,273){\circle*{10}} \put(146,238){\line(1,1){59}}
\put(131,249){\line(1,1){59}} \put(205,296){\line(0,1){40}}
\put(131,210){\line(0,1){40}} \put(242,210){\line(0,1){40}}
\put(316,297){\line(0,1){40}} \put(280,273){\circle*{10}}
\put(279,291){\line(0,-1){35}} \put(279,292){\line(-1,0){22}}
\put(302,256){\line(-1,0){22}} \put(316,297){\line(-3,-1){59}}
\put(302,270){\line(-3,-1){59}} \put(257,277){\line(0,1){15}}
\put(302,257){\line(0,1){12}} \put(103,268){\makebox(22,19){=}}
\put(214,268){\makebox(22,19){=}} \put(28,265){\makebox(17,19){a}}
\end{picture}}

\noindent
{\bf Figure 11 - $UD= DU$} \vspace{3mm}

The functor $F: Tang \longrightarrow Cat(A)$ is defined by
replacing each crossing in a tangle $Q$ by a flat crossing that is
decorated with a corresponding Yang-Baxter element as shown in
Figures 12 and 13. The resulting flat diagram is then a morphism
in $Cat(A).$ This serves to define $F(Q).$ Figure 12 shows how $F$
is naturally defined for vertically oriented crossings. Figure 13
shows how the switchback move implies the definitions of $F$ on
horizontally oriented crossings.  Then in Figure 14 we point out
how the two ways of performing the switchback move are compatible
through our axiomatic assumptions that $(U \otimes U)\rho = \rho$
and $(D \otimes D)\rho = \rho$. This explains and
justifies the first axiom for a quantum algebra. \vspace{3mm}

\noindent {\bf Remark.} The reader should note the compatibility of these 
symmetry axioms about 
$D \otimes D$ and $U \otimes U$ with the categorical turn axiom that states (as in Figure 12)
that the downward versions of the braiding operators (images under $F$ of the crossings)
are exactly the 180 degree turns of the upward versions. We leave as an exercise in 
bead sliding to show that the algebraic symmetry axioms plus the
categorical bead-sliding axioms imply the categorical turn axiom. 
\vspace{3mm}  

 {\tt    \setlength{\unitlength}{0.92pt}
\begin{picture}(349,152)
\thinlines    \put(238,46){\line(0,-1){16}}
              \put(56,48){\line(0,-1){16}}
              \put(236,121){\line(0,-1){16}}
              \put(53,128){\line(0,-1){16}}
              \put(12,93){\vector(3,4){36}}
              \put(52,94){\vector(-3,4){17}}
              \put(28,124){\vector(-3,4){12}}
              \put(121,109){\circle*{10}}
              \put(106,109){\circle*{10}}
              \put(89,101){\makebox(13,13){e}}
              \put(129,101){\makebox(15,15){e'}}
              \put(53,120){\vector(1,0){38}}
              \put(64,103){\makebox(15,16){F}}
              \put(99,97){\vector(3,4){31}}
              \put(128,98){\vector(-3,4){30}}
              \put(144,17){\vector(-3,4){30}}
              \put(115,16){\vector(3,4){31}}
              \put(68,21){\makebox(15,16){F}}
              \put(56,39){\vector(1,0){38}}
              \put(137,44){\circle*{10}}
              \put(124,44){\circle*{10}}
              \put(50,13){\vector(-3,4){35}}
              \put(13,14){\vector(3,4){16}}
              \put(35,41){\vector(3,4){14}}
              \put(100,35){\makebox(17,17){E}}
              \put(145,37){\makebox(14,15){E'}}
              \put(199,61){\vector(3,-4){35}}
              \put(235,59){\vector(-3,-4){15}}
              \put(215,34){\vector(-3,-4){17}}
              \put(226,138){\vector(-3,-4){38}}
              \put(190,138){\vector(3,-4){15}}
              \put(211,113){\vector(3,-4){20}}
              \put(247,96){\makebox(15,16){F}}
              \put(236,113){\vector(1,0){38}}
              \put(271,118){\makebox(15,15){e'}}
              \put(314,118){\makebox(13,13){e}}
              \put(292,125){\circle*{10}}
              \put(307,124){\circle*{10}}
              \put(292,32){\makebox(14,15){E'}}
              \put(322,30){\makebox(17,17){E}}
              \put(304,28){\circle*{10}}
              \put(321,29){\circle*{10}}
              \put(239,38){\vector(1,0){38}}
              \put(251,20){\makebox(15,16){F}}
              \put(284,135){\vector(3,-4){33}}
              \put(314,135){\vector(-3,-4){31}}
              \put(324,56){\vector(-3,-4){31}}
              \put(298,58){\vector(3,-4){33}}
\end{picture}}

\noindent
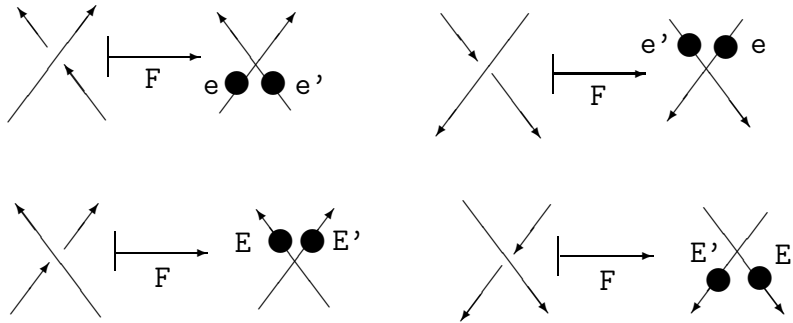
{\bf Figure 12 - The Functor $F$ defined on vertical crossings}
\vspace{3mm}

 {\tt    \setlength{\unitlength}{0.92pt}
\begin{picture}(427,430)
\thinlines    \put(112,21){\vector(0,1){20}}
              \put(111,41){\vector(1,1){15}}
              \put(122,164){\vector(-1,1){15}}
              \put(123,271){\vector(-1,-1){13}}
              \put(347,64){\makebox(18,18){E'}}
              \put(346,18){\makebox(21,21){DE}}
              \put(287,46){\makebox(12,11){=}}
              \put(291,152){\makebox(12,11){=}}
              \put(179,48){\makebox(12,11){=}}
              \put(178,153){\makebox(12,11){=}}
              \put(249,42){\makebox(21,24){E'}}
              \put(210,41){\makebox(18,21){E}}
              \put(247,153){\makebox(17,18){e'}}
              \put(207,155){\makebox(15,14){e}}
              \put(244,58){\circle*{8}}
              \put(231,59){\circle*{8}}
              \put(240,158){\circle*{8}}
              \put(228,159){\circle*{8}}
              \put(223,39){\vector(1,1){30}}
              \put(248,149){\vector(-1,1){31}}
              \put(273,103){\vector(0,-1){62}}
              \put(273,42){\vector(-1,0){24}}
              \put(249,42){\vector(-1,1){28}}
              \put(220,69){\vector(-1,0){17}}
              \put(203,68){\vector(0,-1){56}}
              \put(223,10){\vector(0,1){29}}
              \put(253,68){\vector(0,1){36}}
              \put(274,197){\vector(0,-1){49}}
              \put(273,148){\vector(-1,0){24}}
              \put(216,180){\vector(-1,0){19}}
              \put(197,180){\vector(0,-1){55}}
              \put(221,126){\vector(0,1){27}}
              \put(222,154){\vector(1,1){25}}
              \put(246,178){\vector(0,1){17}}
              \put(146,75){\vector(0,1){36}}
              \put(134,63){\vector(1,1){12}}
              \put(139,175){\vector(0,1){17}}
              \put(114,150){\vector(1,1){25}}
              \put(113,123){\vector(0,1){27}}
              \put(96,75){\vector(0,-1){56}}
              \put(113,76){\vector(-1,0){17}}
              \put(142,49){\vector(-1,1){28}}
              \put(166,49){\vector(-1,0){24}}
              \put(166,110){\vector(0,-1){62}}
              \put(89,176){\vector(0,-1){55}}
              \put(107,178){\vector(-1,0){19}}
              \put(143,142){\vector(-1,1){16}}
              \put(168,142){\vector(-1,0){24}}
              \put(169,191){\vector(0,-1){49}}
              \put(383,278){\makebox(34,32){UE}}
              \put(323,344){\makebox(29,28){Ue'}}
              \put(68,38){\makebox(12,11){=}}
              \put(67,157){\makebox(12,11){=}}
              \put(367,161){\circle*{8}}
              \put(366,142){\circle*{8}}
              \put(348,61){\circle*{8}}
              \put(347,45){\circle*{8}}
              \put(37,49){\vector(-1,1){16}}
              \put(58,27){\vector(-1,1){16}}
              \put(35,159){\vector(-1,-1){14}}
              \put(59,183){\vector(-1,-1){16}}
              \put(374,73){\vector(-1,-1){39}}
              \put(373,35){\vector(-1,1){38}}
              \put(376,171){\vector(-1,-1){39}}
              \put(375,133){\vector(-1,1){38}}
              \put(59,66){\vector(-1,-1){39}}
              \put(56,146){\vector(-1,1){38}}
              \put(102,358){\vector(0,-1){19}}
              \put(136,385){\vector(-4,-3){34}}
              \put(136,417){\vector(0,-1){31}}
              \put(157,355){\vector(0,1){65}}
              \put(134,354){\vector(1,0){22}}
              \put(119,368){\vector(1,-1){14}}
              \put(100,387){\vector(1,-1){15}}
              \put(75,387){\vector(1,0){24}}
              \put(75,343){\vector(0,1){43}}
              \put(12,254){\vector(3,4){39}}
              \put(16,307){\vector(3,-4){15}}
              \put(37,280){\vector(3,-4){20}}
              \put(33,380){\vector(3,4){16}}
              \put(13,355){\vector(3,4){14}}
              \put(12,401){\vector(3,-4){35}}
              \put(218,342){\vector(0,1){43}}
              \put(219,386){\vector(1,0){24}}
              \put(277,353){\vector(1,0){22}}
              \put(299,354){\vector(0,1){65}}
              \put(276,415){\vector(0,-1){31}}
              \put(276,383){\vector(-4,-3){34}}
              \put(242,356){\vector(0,-1){19}}
              \put(244,385){\vector(1,-1){32}}
              \put(269,377){\circle*{8}}
              \put(253,377){\circle*{8}}
              \put(236,369){\makebox(11,12){e'}}
              \put(277,369){\makebox(12,13){e}}
              \put(343,338){\vector(3,4){51}}
              \put(341,405){\vector(3,-4){52}}
              \put(359,358){\circle*{8}}
              \put(357,385){\circle*{8}}
              \put(334,378){\makebox(14,13){e}}
              \put(51,371){\makebox(12,11){=}}
              \put(308,373){\makebox(12,11){=}}
              \put(167,378){\vector(1,0){38}}
              \put(177,358){\makebox(16,16){F}}
              \put(187,263){\makebox(16,16){F}}
              \put(177,283){\vector(1,0){38}}
              \put(317,278){\makebox(12,11){=}}
              \put(378,265){\circle*{8}}
              \put(377,293){\circle*{8}}
              \put(340,312){\vector(3,-4){52}}
              \put(341,246){\vector(3,4){51}}
              \put(276,268){\circle*{8}}
              \put(264,269){\circle*{8}}
              \put(254,290){\vector(1,-1){32}}
              \put(252,261){\vector(0,-1){19}}
              \put(287,288){\vector(-4,-3){34}}
              \put(286,320){\vector(0,-1){31}}
              \put(309,259){\vector(0,1){65}}
              \put(287,258){\vector(1,0){22}}
              \put(229,291){\vector(1,0){24}}
              \put(228,247){\vector(0,1){43}}
              \put(85,248){\vector(0,1){43}}
              \put(85,292){\vector(1,0){24}}
              \put(143,259){\vector(1,0){22}}
              \put(166,260){\vector(0,1){65}}
              \put(146,323){\vector(0,-1){31}}
              \put(109,257){\vector(0,-1){19}}
              \put(109,290){\vector(1,-1){32}}
              \put(145,292){\vector(-1,-1){16}}
              \put(57,272){\makebox(13,13){=}}
              \put(240,262){\makebox(16,15){E'}}
              \put(281,265){\makebox(15,15){E}}
              \put(384,257){\makebox(16,15){E'}}
              \put(377,150){\makebox(17,21){De'}}
              \put(372,128){\makebox(25,26){e}}
              \put(166,388){\line(0,-1){20}}
              \put(177,292){\line(0,-1){20}}
\end{picture}}

\noindent
{\bf Figure 13 - Horizontal Crossings} \vspace{3mm}

 {\tt    \setlength{\unitlength}{0.92pt}
\begin{picture}(407,267)
\thinlines    \put(302,133){$U^{-1}e$}
              \put(137,138){\line(0,-1){22}}
              \put(165,226){\line(0,-1){22}}
              \put(27,120){\makebox(12,11){=}}
              \put(278,111){\makebox(12,11){=}}
              \put(360,135){\circle*{8}}
              \put(360,108){\circle*{8}}
              \put(343,155){\vector(3,-4){52}}
              \put(344,88){\vector(3,4){51}}
              \put(148,107){\makebox(16,16){F}}
              \put(138,127){\vector(1,0){38}}
              \put(243,109){\makebox(13,13){e'}}
              \put(208,111){\makebox(11,11){e}}
              \put(237,112){\circle*{8}}
              \put(225,113){\circle*{8}}
              \put(242,106){\vector(-1,1){28}}
              \put(190,153){\vector(0,-1){47}}
              \put(192,107){\vector(1,0){29}}
              \put(220,107){\vector(2,3){16}}
              \put(237,132){\vector(1,0){26}}
              \put(262,132){\vector(0,-1){54}}
              \put(243,79){\vector(0,1){26}}
              \put(215,135){\vector(0,1){18}}
              \put(69,131){\vector(0,1){18}}
              \put(83,117){\vector(-1,1){12}}
              \put(100,101){\vector(-1,1){12}}
              \put(100,75){\vector(0,1){26}}
              \put(120,127){\vector(0,-1){54}}
              \put(94,128){\vector(1,0){26}}
              \put(78,104){\vector(2,3){16}}
              \put(49,103){\vector(1,0){29}}
              \put(48,150){\vector(0,-1){47}}
              \put(101,194){\vector(0,-1){19}}
              \put(136,222){\vector(-4,-3){34}}
              \put(135,255){\vector(0,-1){31}}
              \put(156,192){\vector(0,1){65}}
              \put(134,191){\vector(1,0){22}}
              \put(118,205){\vector(1,-1){14}}
              \put(99,224){\vector(1,-1){15}}
              \put(74,224){\vector(1,0){24}}
              \put(74,179){\vector(0,1){43}}
              \put(32,217){\vector(3,4){16}}
              \put(12,192){\vector(3,4){14}}
              \put(11,238){\vector(3,-4){35}}
              \put(217,179){\vector(0,1){43}}
              \put(218,223){\vector(1,0){24}}
              \put(276,190){\vector(1,0){22}}
              \put(298,191){\vector(0,1){65}}
              \put(275,252){\vector(0,-1){31}}
              \put(275,220){\vector(-4,-3){34}}
              \put(241,193){\vector(0,-1){19}}
              \put(243,222){\vector(1,-1){32}}
              \put(268,214){\circle*{8}}
              \put(252,214){\circle*{8}}
              \put(235,206){\makebox(11,12){e'}}
              \put(276,206){\makebox(12,13){e}}
              \put(342,175){\vector(3,4){51}}
              \put(341,242){\vector(3,-4){52}}
              \put(358,195){\circle*{8}}
              \put(358,222){\circle*{8}}
              \put(50,208){\makebox(12,11){=}}
              \put(307,210){\makebox(12,11){=}}
              \put(166,215){\vector(1,0){38}}
              \put(176,195){\makebox(16,16){F}}
              \put(114,10){\makebox(211,47){$(U \otimes U) \rho = \rho$}}
              \put(328,220){$e$}
              \put(315,194){$Ue'$}
              \put(323,105){$e'$}
\end{picture}}

\noindent
{\bf Figure 14 - The relations $(U \otimes U)\rho = \rho$ and $(D
\otimes D)\rho = \rho$} \vspace{3mm}

In order to see that $F$ is well-defined it remains to verify
invariance under the reverse Reidemeister two move. In Figure 15
we show that this invariance is equivalent to the equation

$$1_{A} \otimes 1_{A} = eDE \otimes E'Ue'$$

(summation convention of the pairs $E,E'$ and $e,e'$). This, in
turn, is equivalent to the condition that

$$[(1_{A} \otimes U) \rho][(D \otimes 1_{A^{op}}) \rho^{-1}] =
1_{A \otimes A^{op}}.$$

\noindent This proves that the functor $F$ is well-defined and it
gives us the motivation for the second axiom for an oriented
quantum algebra. \vspace{3mm}

 {\tt    \setlength{\unitlength}{0.92pt}
\begin{picture}(377,270)
\thinlines    \put(145,210){$F$}
              \put(140,218){\line(0,-1){27}}
              \put(141,204){\vector(1,0){40}}
              \put(212,246){\circle*{8}}
              \put(211,215){\circle*{8}}
              \put(211,182){\circle*{8}}
              \put(214,162){\circle*{8}}
              \put(90,256){\vector(3,-4){40}}
              \put(83,201){\vector(3,4){18}}
              \put(113,236){\vector(3,4){17}}
              \put(130,202){\vector(-3,-4){39}}
              \put(130,151){\vector(-1,1){19}}
              \put(105,176){\vector(-1,1){23}}
              \put(243,202){\vector(-3,-4){39}}
              \put(203,256){\vector(3,-4){40}}
              \put(195,200){\vector(1,1){53}}
              \put(240,154){\vector(-1,1){45}}
              \put(12,255){\vector(0,-1){103}}
              \put(42,153){\vector(0,1){101}}
              \put(53,195){\makebox(12,14){=}}
              \put(336,151){\vector(0,1){101}}
              \put(306,253){\vector(0,-1){103}}
              \put(337,187){\circle*{8}}
              \put(337,229){\circle*{8}}
              \put(306,230){\circle*{8}}
              \put(306,188){\circle*{8}}
\thicklines   \put(53,90){\vector(1,0){45}}
\thinlines    \put(164,23){\vector(0,1){101}}
              \put(132,22){\vector(0,1){102}}
              \put(214,22){\vector(0,1){102}}
              \put(245,24){\vector(0,1){101}}
              \put(215,97){\circle*{8}}
              \put(245,97){\circle*{8}}
              \put(214,61){\circle*{8}}
              \put(246,61){\circle*{8}}
              \put(286,222){\makebox(13,14){e}}
              \put(344,176){\makebox(21,20){E'}}
              \put(191,238){\makebox(13,14){e}}
              \put(184,172){\makebox(21,20){E'}}
              \put(257,195){\makebox(12,14){=}}
              \put(250,51){\makebox(21,20){E'}}
              \put(194,55){\makebox(13,14){e}}
              \put(176,71){\makebox(12,14){=}}
\thicklines   \put(114,10){\framebox(165,123){}}
\thinlines    \put(277,177){\makebox(24,24){DE}}
              \put(179,150){\makebox(24,24){DE}}
              \put(183,208){\makebox(22,23){Ue'}}
              \put(345,216){\makebox(22,23){Ue'}}
              \put(186,82){\makebox(24,24){DE}}
              \put(250,85){\makebox(22,23){Ue'}}
\end{picture}}

\noindent
{\bf Figure 15 - The Reverse Type Two Move} \vspace{3mm}

\subsection{Bead Sliding} Morphisms in the category $Cat(A)$ are
precisely the morphisms in $Flat$ decorated by elements of the
algebra $A$. The rules of interaction of cups and caps with
algebra elements amount to this: If an algebra element $a$
decorates the right side of a a cap then it can be moved to the
left side of that cap at the expense of relabeling it as $U(a)$
if the left side of the cap is a rising orientation, and $D(a)$ if
the left side of the cap is a falling orientation. If an algebra
element $a$ decorates the left side of a a cup then it can be
moved to the right side of that cup at the expense of relabeling
it as $U(a)$ if the right side of the cup is a rising orientation,
and $D(a)$ if the right side of the cup is a falling orientation.
\vspace{3mm}

This amounts to a counterclockwise turn or {\em bead slide} taking
$a$ around the cap and changing it into $U(a)$ or $D(a)$. This
rule of transformation is uniform. Clockwise turns correspond to
applications of the automorphisms $U^{-1}$ and $D^{-1}$, while
counterclockwise turns correspond to applications of $U$ and $D$.
The upshot of these transformations is that, given a morphism in
$Cat(A)$, we can move all the decorations on a given component to
any single vertical segment of that component. In particular, this
means that in the case of a $1-1$ tangle, we can move all the
algebra to the top of the tangle. What remains below the algebra
is an immersed curve that can be regularly homotoped to a string
of curls (as illustrated in Figure 11). Each curl is a special
morphism in this category, not neccessarily corresponding to an
element of the algebra. It is convenient to label the curls $G$
and $G^{-1}$ as shown in Figure 11. Then we get an algebraic
expression for every morphism in the form $G^{n}w$ where $w$ is an
expression in the algebra $A$.  Note that in general $w$ is a
summation of products since the decoration of each crossing
consists in a sum in $A \otimes A$. \vspace{3mm}

\noindent For a $1-1$ tangle $T$ we let $Inv(T)$ denote the
expression $G^{n}w$. Up to equivalence in the oriented algebra A
(possibly augmented by an element corresponding to $G$) the
expression $Inv(T)$ is an invariant of the regular isotopy class of
the tangle. The following proposition shows what happens when we 
reverse the orientation of the tangle T. \vspace{3mm}

\noindent {\bf Proposition.} Let $T$ be a $1-1$ tangle of one component 
(a knot with ends). Let $r(T)$ denote the result of reversing the
orientation of $T$. If $Inv(T) = G^{n}w$ then $Inv(r(T)) = G^{-n}r(w)$
where $r(w)$ is obtained by reversing the order of
the products in $w$.
\vspace{3mm}

\noindent {\bf Proof.} First note that we can arrange the diagram for the
tangle so that all the crossings are vertical. It then follows
from our conventions that reversing the orientation of the diagram
does not affect the decoration of the crossings with algebra. See
Figure 12 for an illustration.  The cup and cap operators are
obtained for the reverse orientation by interchanging $U$ and $D$.
Note that if (for the sake of argument) the pair of signifiers for
the $\rho$ at a crossing are $E$ and $E'$ then, when the beads are
slid to the top of the tangle, $E$ and $E'$ will differ by a
composition of the operators $U$ and $D$ that is obtained from
going around a closed loop (from down to down or from up to up).
Such a composition of operators is neccessarily a power of $t = UD
= DU$. We know that $(U \otimes U) \rho = \rho$  and that $(D
\otimes D) \rho^{-1} = \rho^{-1}.$ Note that $U \otimes U$ is an
automorphism of $A \otimes A$. Hence 
$ \rho^{-1} = ((U \otimes U)\rho)^{-1} = (U \otimes U)\rho^{-1}.$
Similarly, $(D \otimes D)\rho = \rho.$ The upshot of this 
remark is that in the algebra expression $Inv(T)$ if $E$ and $E'$
receive identical operators as compositions of $U$ and $D$, then 
these operators can be removed from both of them. Since, referring
specifically to $E$ and $E'$ in the discussion above, we have
that $E$ and $E'$ will differ by a power of $t = DU = UD$, it 
follows that the entire algebra expression $Inv(T)$ involves only 
the operator $t$ (after removing common operators from the pairs
$E$, $E'$ and $e$, $e'$. The Lemma follows easily from these 
observations. //
\vspace{3mm}

\noindent {\bf Corollary.} The knot invariant derived via $1-1$ tangles from
an oriented quantum algebra is equivalent to an invariant 
derived from a standard oriented quantum algebra obtained by
replacing $U$ by $t = UD = DU$ and replacing $D$ by the identity
automorphism.
\vspace{3mm}

\noindent
{\bf Proof.} It is not hard to see that if $(A,\rho,D,U)$ is an 
oriented quantum algebra, then $(A,\rho,1,DU)$ is also a quantum 
algebra. The method of the proof of the proposition shows that 
these two algebras yield the same invariants. //
\vspace{3mm}

 {\tt    \setlength{\unitlength}{0.92pt}
\begin{picture}(416,485)
\thinlines    \put(55,146){\line(0,-1){26}}
              \put(158,374){\line(0,-1){27}}
              \put(329,348){\vector(3,4){67}}
              \put(284,436){\vector(0,-1){172}}
              \put(332,428){\vector(-4,-3){70}}
              \put(329,283){\vector(0,1){67}}
              \put(307,282){\vector(1,0){21}}
              \put(262,318){\vector(4,-3){46}}
              \put(286,264){\vector(1,0){82}}
              \put(263,375){\vector(0,-1){58}}
              \put(368,264){\vector(0,1){167}}
              \put(367,429){\vector(-1,0){35}}
              \put(245,437){\vector(1,0){40}}
              \put(13,434){\vector(1,0){40}}
              \put(140,405){\vector(1,1){25}}
              \put(135,425){\vector(-1,0){35}}
              \put(136,261){\vector(0,1){167}}
              \put(50,390){\vector(-1,-1){18}}
              \put(31,372){\vector(0,-1){58}}
              \put(54,262){\vector(1,0){82}}
              \put(30,315){\vector(4,-3){46}}
              \put(75,279){\vector(1,0){21}}
              \put(97,280){\vector(0,1){67}}
              \put(99,425){\vector(-4,-3){39}}
              \put(54,434){\vector(0,-1){132}}
              \put(53,291){\vector(0,-1){29}}
              \put(97,347){\vector(2,3){34}}
\thicklines   \put(159,361){\vector(1,0){60}}
              \put(249,184){\circle*{8}}
              \put(248,207){\circle*{8}}
              \put(250,143){\circle*{8}}
              \put(249,164){\circle*{8}}
              \put(249,100){\circle*{8}}
              \put(249,122){\circle*{8}}
              \put(380,414){\circle*{8}}
              \put(368,414){\circle*{8}}
              \put(298,291){\circle*{8}}
              \put(284,291){\circle*{8}}
              \put(283,381){\circle*{8}}
              \put(273,382){\circle*{8}}
              \put(347,405){\makebox(15,17){E}}
              \put(390,406){\makebox(16,16){E'}}
              \put(252,384){\makebox(15,15){F'}}
              \put(289,371){\makebox(17,18){F}}
              \put(261,285){\makebox(15,16){K'}}
              \put(304,286){\makebox(17,19){K}}
              \put(55,134){\vector(1,0){57}}
\thinlines    \put(211,10){\vector(0,1){92}}
              \put(210,102){\vector(-1,0){34}}
              \put(176,102){\vector(0,-1){41}}
              \put(178,61){\vector(1,0){71}}
              \put(248,61){\vector(0,1){164}}
              \put(258,199){\makebox(12,19){E'}}
              \put(257,174){\makebox(22,18){UK}}
              \put(259,157){\makebox(22,15){UF'}}
              \put(259,136){\makebox(23,17){UDE}}
              \put(260,113){\makebox(27,20){UDUK'}}
              \put(257,91){\makebox(28,20){UDUF}}
              \put(177,326){\makebox(20,26){F}}
              \put(56,140){\makebox(54,28){Bead Slide}}
              \put(165,430){\vector(0,1){37}}
              \put(397,438){\vector(0,1){37}}
              \put(246,246){\vector(0,1){191}}
              \put(12,243){\vector(0,1){191}}
\end{picture}}

\noindent
{\bf Figure 16 - Bead Sliding Trefoil}
\vspace{3mm}

\noindent
{\bf Remark.} In Figure 16 we have 

$$w(T) = G^{-1}(UDUF)(UDUK')(UDE)(UF')(UK)(E')$$ 
$$= G^{-1}(DUF)(DUK')(UDE)F'KE' = G^{-1}(tF)(tK')(tE)F'KE'.$$

\noindent This example is a concise illustration of the content of the
Proposition. This proposition suggests the conjecture that the $1-1$ tangle
invariants can detect the difference between some 
knots and their reversals.
\vspace{3mm}

\section{Matrix Models} This section will show how matrix
representations of a quantum algebra (or a matrix quantum algebra)
give rise to invariants that can be construed directly as state
summations via the vertex weights from cup, cap and crossing
matrices. Thus we get a double description in terms of bead
sliding and in terms of the state sums. We will {\em also} discuss
how our theory of quantum algebras corresponds to the usual way of
augmenting a solution to the Yang-Baxter equation to produce a
quantum link invariant.  This is a case of taking an oriented
quantum algebra associated with a solution to the Yang-Baxter
equation. We shall restrict the discussion to balanced oriented
quantum algebras where $U=D=T$.  We modify our discussion
for the general case and for other specific cases of quantum
algebras in \cite{KRNexTNexT}. \vspace{3mm}

In a matrix model for a knot or link invariant we are given a
diagram for the knot or link $K$ that is arranged with respect to
a vertical direction so that there are a finite number of
transverse directions to the vertical that have critical points
(maxima, minima or crossings) and these critical points are
separated. Thus $K$ is given as a morphism in the tangle category.
We then traverse the diagram for $K$, marking one point on each
arc in the diagram from critical point to critical point. This
divides the diagram (by deleting the marked points) into a
collection of generators of the tangle category (cups, caps and
crossings). See Figure 17. \vspace{3mm}

 {\tt    \setlength{\unitlength}{0.92pt}
\begin{picture}(282,312)
\thinlines    \put(257,259){\vector(-1,1){12}}
              \put(130,270){\vector(-1,-1){16}}
              \put(154,292){\vector(-1,-1){16}}
              \put(48,283){\vector(1,1){13}}
              \put(255,286){\framebox(13,14){j}}
              \put(100,240){\framebox(13,14){j}}
              \put(154,70){\framebox(13,14){j}}
              \put(63,252){\framebox(13,14){i}}
              \put(214,288){\framebox(13,14){i}}
              \put(128,70){\framebox(13,14){i}}
              \put(162,64){\circle*{8}}
              \put(124,65){\circle*{8}}
              \put(155,240){\framebox(10,13){h}}
              \put(254,216){\framebox(10,13){h}}
              \put(259,243){\framebox(13,15){g}}
              \put(228,215){\framebox(13,15){g}}
              \put(12,252){\framebox(12,13){e}}
              \put(216,241){\framebox(12,15){f}}
              \put(175,214){\framebox(12,15){f}}
              \put(153,217){\framebox(12,13){e}}
              \put(156,287){\framebox(12,13){d}}
              \put(36,193){\framebox(12,13){d}}
              \put(106,288){\framebox(11,12){c}}
              \put(115,192){\framebox(11,12){c}}
              \put(92,193){\framebox(10,11){b}}
              \put(64,290){\framebox(10,11){b}}
              \put(10,288){\framebox(12,11){a}}
              \put(11,194){\framebox(12,11){a}}
              \put(117,290){\vector(1,-1){36}}
              \put(24,262){\vector(1,1){13}}
              \put(62,265){\vector(-4,3){37}}
              \put(258,216){\vector(0,-1){27}}
              \put(258,189){\vector(-1,0){25}}
              \put(233,191){\vector(0,1){26}}
              \put(95,206){\vector(0,1){26}}
              \put(95,233){\vector(1,0){25}}
              \put(121,232){\vector(0,-1){27}}
              \put(240,275){\vector(-1,1){12}}
              \put(254,285){\vector(-1,-1){27}}
              \put(156,193){\vector(0,1){26}}
              \put(181,191){\vector(-1,0){25}}
              \put(181,218){\vector(0,-1){27}}
              \put(42,234){\vector(0,-1){27}}
              \put(16,235){\vector(1,0){25}}
              \put(16,208){\vector(0,1){26}}
              \put(232,25){\framebox(10,13){h}}
              \put(177,26){\framebox(13,15){g}}
              \put(114,24){\framebox(12,15){f}}
              \put(194,118){\framebox(12,13){d}}
              \put(55,25){\framebox(12,13){e}}
              \put(169,118){\framebox(11,12){c}}
              \put(82,108){\framebox(10,11){b}}
              \put(54,118){\framebox(12,11){a}}
              \put(170,112){\circle*{8}}
              \put(46,123){\circle*{8}}
              \put(201,113){\circle*{8}}
              \put(135,29){\circle*{8}}
              \put(49,31){\circle*{8}}
              \put(168,30){\circle*{8}}
              \put(251,32){\circle*{8}}
              \put(99,110){\circle*{8}}
              \put(92,101){\vector(2,3){17}}
              \put(109,126){\vector(1,0){50}}
              \put(158,125){\vector(1,-1){109}}
              \put(265,17){\vector(-1,0){79}}
              \put(141,53){\vector(-4,3){91}}
              \put(44,124){\vector(0,1){38}}
              \put(44,162){\vector(1,0){194}}
              \put(237,161){\vector(-3,-4){46}}
              \put(181,91){\vector(-3,-4){58}}
              \put(122,12){\vector(-1,0){87}}
              \put(36,12){\vector(2,3){45}}
              \put(185,17){\vector(-4,3){33}}
\end{picture}}

\noindent
{\bf Figure 17 -- A knot $K$ with marked points} \vspace{3mm}

\noindent We further assume that a matrices with  entries in an
appropriate commutative ring $k$ have been assigned to each of the
different orientation types of cups, caps and crossings. With such
an assigment of matrices, we can create an evaluation $Z(K)$ of a
given marked diagram by the following algorithm: Let $I$ denote
the index set for the individual indices on the cap, cap and
crossing matrices. (Cups and caps have two indices while crossings
have four indices.) A {\em coloring} of a marked diagram is an
arbitrary assignment of indices to the marked points on the
diagram. A colored diagram then has indices assigned to each of
its component cup, cap and crossing matrices. {\em We define Z(K)
to be the sum over all colorings of the products of these matrix
entries.} Thus in the diagram shown in Figure 17 we have that

$$Z(K) = \Sigma_{a,...,h} [\overline{M}_{ad}^{>}] [
\overline{S}^{cd}_{ef}] [\underline{S}^{cd}_{gh}] [|R^{fg}_{kl}]
[\underline{M}^{ck}_{>}] [\underline{M}^{eh}_{<}].$$

\noindent Here $[\overline{M}^{>}]$ stands for a right-oriented
(clockwise) cap, $[\overline{M}^{<}]$ stands for a left-oriented
(counter-clockwise) cap, $[\underline{M}_{>}]$  stands for a
counter-clockwise cup, $[\underline{M}_{<}]$ stands for a
clockwise cup. $\overline{R}$ is the positive upward pointing
crossing matrix and $\overline{S}$ is its inverse. A bar below the
$R$ connotes a downward-pointing crossing and a bar to the left or
to the right connotes a crossing that points to the left or to the
right respectively. Thus we have the following list of matrices
that are directly associated with the link diagram:

$$\overline{M}^{>},\overline{M}^{<}, \underline{M}_{>}, 
\underline{M}_{<}, \overline{R}, \underline{R}, |R,
R|, \overline{S}, \underline{S}, |S, S|$$
\vspace{3mm}

 {\tt    \setlength{\unitlength}{0.92pt}
\begin{picture}(397,350)
\thinlines    \put(269,81){$\underline{S}$}
              \put(67,86){$\overline{S}$}
              \put(256,203){$\underline{R}$}
              \put(58,207){$\overline{R}$}
              \put(55,257){$\underline{M}_{>}$}
              \put(254,256){$\underline{M}_{<}$}
              \put(57,316){$\overline{M}^{>}$}
              \put(284,23){$S|$}
              \put(274,149){$R|$}
              \put(77,29){$|S$}
              \put(69,155){$|R$}
              \put(39,33){\vector(-3,2){15}}
              \put(61,17){\vector(-3,2){16}}
              \put(58,46){\vector(-1,-1){30}}
              \put(30,154){\vector(-1,-1){14}}
              \put(50,173){\vector(-1,-1){12}}
              \put(49,146){\vector(-4,3){31}}
              \put(239,24){\vector(1,-1){13}}
              \put(221,40){\vector(1,-1){11}}
              \put(221,12){\vector(1,1){29}}
              \put(220,99){\vector(1,-1){30}}
              \put(229,208){\vector(1,-1){11}}
              \put(210,225){\vector(1,-1){11}}
              \put(240,224){\vector(-1,-1){26}}
              \put(42,93){\vector(1,1){12}}
              \put(21,73){\vector(1,1){13}}
              \put(51,75){\vector(-1,1){30}}
              \put(22,218){\vector(-1,1){12}}
              \put(45,198){\vector(-1,1){14}}
              \put(15,198){\vector(3,4){25}}
              \put(14,318){\vector(0,1){19}}
              \put(13,338){\vector(1,0){21}}
              \put(34,339){\vector(0,-1){21}}
              \put(234,313){\vector(0,1){19}}
              \put(233,332){\vector(-1,0){17}}
              \put(217,331){\vector(0,-1){19}}
              \put(14,274){\vector(0,-1){19}}
              \put(16,255){\vector(1,0){19}}
              \put(35,257){\vector(0,1){18}}
              \put(236,269){\vector(0,-1){20}}
              \put(234,249){\vector(-1,0){18}}
              \put(216,248){\vector(0,1){20}}
              \put(230,82){\vector(-1,-1){14}}
              \put(253,103){\vector(-1,-1){15}}
              \put(212,171){\vector(1,-1){35}}
              \put(210,134){\vector(1,1){18}}
              \put(233,157){\vector(1,1){15}}
              \put(254,312){$\overline{M}^{<}$}
\end{picture}}

\noindent
{\bf Figure 18 - Matrix Notations}
\vspace{3mm}

In the models all these oriented "matrices" (as above and in Figure 18) will be
defined by a smaller collection of ordinary
matrices with standard multiplication convention. We shall call
these oriented matrices of Figure 18 the {\em diagram matrices} since they can
be read directly from a decorated Morse link diagram in the
process of translating from topology to algebra. We shall call the
smaller collection of standard matrices the {\em background matrices} for
the matrix model. There will be a single background matrix $M$
(written with lower indices $M_{ab}$) and its inverse $M^{-1}$ so
that $\Sigma_{i} M_{ai}M^{-1}_{ib} = \delta_{ab}.$ The background
matrices $R$ and $S$ written with indices in the form
$R^{ab}_{cd}$ and $S^{ab}_{cd}$, and $\Sigma_{ij}
R^{ij}_{cd}S^{ab}_{ij} = \delta^{a}_{c} \delta^{b}_{d}.$ Thus, if
we think of $ab$ as a single index and write $R^{ab}_{cd} = R_{ab,
cd}$ and $S^{ab}_{cd}$ then $R$ and $S$ multiply in the standard
matrix convention where the left lower index is the input index
and the right lower index is the output index. In this convention
we can write $RS=SR=I$ where $I$ denotes the identity matrix of
this dimension. \vspace{3mm}

The rewrite definitions of the diagram matrices in terms of the
background matrices are as follows:

$$\overline{M}^{>}_{ab} = M_{ab}$$ $$\overline{M}^{<}_{ab}=
M^{-1}_{ba}$$ $$\underline{M}_{>}^{ab}= M^{-1}_{ab}$$
$$\underline{M}_{<}^{ab} = M_{ba}$$ $$\overline{R}^{ab}_{cd} =
R^{ab}_{cd}$$ $$\underline{R}^{ab}_{cd} = R^{dc}_{ba}$$
$$\overline{S}^{ab}_{cd} = S^{ab}_{cd}$$ $$\underline{S}^{ab}_{cd}
= S^{dc}_{ba}$$ $$|R^{ab}_{cd} = \overline{M}^{<}_{ci}
\overline{R}^{ia}_{dj} \underline{M}_{<}^{jb} =
\underline{M}_{<}^{jb} \underline{R}^{bj}_{ic}
\overline{M}^{<}_{jd}$$ $$R|^{ab}_{cd} = \overline{M}^{>}_{ci}
\underline{R}^{ia}_{dj} \underline{M}_{>}^{jb} =
\underline{M}_{>}^{jb} \overline{R}^{bj}_{ic}
\overline{M}^{>}_{jd}$$ $$|S^{ab}_{cd} = \overline{M}^{<}_{ci}
\overline{S}^{ia}_{dj} \underline{M}_{<}^{jb} =
\underline{M}_{<}^{jb} \underline{S}^{bj}_{ic}
\overline{M}^{<}_{jd}$$ $$S|^{ab}_{cd} = \overline{M}^{>}_{ci}
\underline{S}^{ia}_{dj} \underline{M}_{>}^{jb} =
\underline{M}_{>}^{jb} \overline{S}^{bj}_{ic}
\overline{M}^{>}_{jd}$$

The last four equations expressing the horizontal versions of the
braiding matrices in terms of the vertical ones and the cup and
cap matrices follow from the switchback move as illustrated in
Figures 5 and 13.  The conditions for invariance under regular isotopy in
the tangle category (Figure 5) translate into conditions on these matrices.
For example, the type three move translates to the Yang-Baxter
Equation in braiding form. Assuming that the matrices satisfy these
conditions, it follows that $Z(K)$ is a regular isotopy invariant
of knots and links. \vspace{3mm}

In fact the same method of assignment gives a functor from the
tangle category to the tensor category associated with the basic
representation module associated with these matrices or
equivalently  to the category of an oriented quantum algebra
$M_{n}(k)$ associated with the $n \times n$ matrices over the ring
$k$ where $n$ is the size of the index set. In working with $n
\times n$ matrices $A$ it is convenient for diagram purposes to
write $A_{i}^{j} = A_{ij}$ where the second half of this equation
denotes the standard convention for matrix indices ($ij$ stands
for row-$i$ and column-$j$) so that $$(AB)_{ij} = \Sigma_{k}
A_{ik}B_{kj}.$$ In writing $A_{i}^{j}$ we indicate that for the
standard upward orientation with repect to the vertical, the input
index for the matrix is at the bottom and the output index for the
matrix is at the top. With these conventions, Figure 10 (interpreted for matrices)
shows that
the automorpism $T:M_{n}(k) \longrightarrow A(n)$ is given by the
equation $$T(A) = MAM^{-1}.$$

\noindent To see this note that $$T(A)_{ij} = T(A)_{i}^{j} =
\overline{M}^{>}_{ik} A^{l}_{k} \underline{M}_{>}^{lj} =
M_{ik}A_{kl}M^{-1}_{lj} = (MAM^{-1})_{ij}.$$

\noindent In this picture the generators of the algebra are the
elementary matrices $E_{a}^{b}$ (with entry $1$ in the a-th row,
b-th column place and zero elsewhere). Here we think of the lower
index on the elementary matrix as the index corresponding to the
entrance to the arrow for the corresponding morphism in
$Cat(M_{n}(k))$ and the upper index corresponds to the exit from the
arrow. In this convention we have $E_{a}^{b}E_{c}^{d} =
\delta^{b}_{c}E_{a}^{d}$ where  $\delta^{b}_{c}$ is the Kronecker
delta (equal to one when $b=c$ and $0$ otherwise). These
conventions are important for specific calculations of these
quantum algebras associated to the matrix models. \vspace{3mm}

In order to complete the relationship between matrix models and our description
of oriented quantum algebras we note that given $\rho \in M_{n}(k) \otimes M_{n}(k)$
We obtain a braiding matrices $R$ and $\overline{R}$ by permuting the upper and lower 
indices of $\rho$ and $\rho^{-1}$ respectively. That is

$$R^{ab}_{cd} = \rho^{ba}_{cd}$$

\noindent and

$$\overline{R}^{ab}_{cd} = (\rho^{-1})^{ab}_{dc}.$$

\noindent These assignments follow directly from the definition of the functor from 
the tangle category to the category of the algebra $M_{n}(k).$  Figure 19 illustrates
the corresponding diagrams. In these diagrams we denote the matrix representation of 
$\rho = e \otimes e'$ by 

$$\rho^{ab}_{cd} = e^{a}_{c} e'^{b}_{d}.$$

\noindent The functor places the signifiers $e$ and $e'$ of $\rho$ on the lines of a 
flat crossing and the resulting diagrammatic matrix is given by 

$$R^{ab}_{cd} = e'^{a}_{d}e^{b}_{c} = e^{b}_{c}e'^{a}_{d}= \rho^{ba}_{cd}$$

\noindent Note that this identity is dependent on the fact that the matrix elements that 
correspond to the signifiers of $\rho$ are commuting scalars in the base ring $k$.
\vspace{3mm}

{\tt    \setlength{\unitlength}{0.92pt}
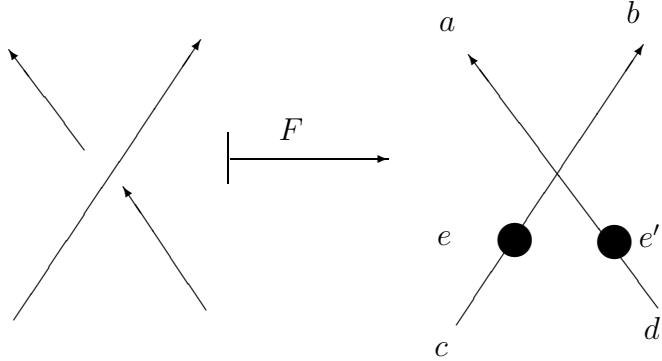
\begin{picture}(308,209)
\thinlines    \put(27,13){$R^{ab}_{cd} = e'^{a}_{d}e^{b}_{c}=\rho^{ba}_{cd}$}
              \put(122,143){$F$}
              \put(101,146){\line(0,-1){21}}
              \put(102,135){\vector(1,0){65}}
              \put(272,61){$d$}
              \put(186,54){$c$}
              \put(265,191){$b$}
              \put(188,188){$a$}
              \put(270,99){$e'$}
              \put(187,100){$e$}
              \put(260,101){\circle*{14}}
              \put(-78,195){\circle*{0}}
              \put(219,102){\circle*{14}}
              \put(278,74){\vector(-3,4){78}}
              \put(195,67){\vector(2,3){77}}
              \put(42,139){\vector(-3,4){31}}
              \put(92,73){\vector(-2,3){34}}
              \put(13,69){\vector(2,3){77}}
\end{picture}}

\noindent
{\bf Figure 19 - Braiding Matrix $R$ and Algebraic Matrix $\rho$}
\vspace{3mm}

\subsection {Matrices and Bead Sliding.} Note that when we use a matrix
model based on a representation of an oriented (balanced) quantum algebra we require
not only a matrix representation of the algebra, but also a matrix representation of the
basic automorphism $T$ of $A$. This data entails the matrices $e^{a}_{b}$ and $e'^{a}_{b}$
that define the matrix representation of $\rho$ and the matrices $E^{a}_{b}$ and 
$E'^{a}_{b}$ that define the matrix representation of $\rho^{-1}$ and the background matrices
$M$ and $M^{-1}$ that define the representation of $T$ via $T(v) = MvM^{-1}.$ Note that we
apply $T$ functionally on the left. It has been our convention to multiply
algebra in the order of its appearance on the oriented lines of the diagram. With our index
conventions for matrices this corresponds to the left-right order of matrix multiplication.
This means that if we regard the functor $Cat(M_{n}(k))$ as containing morphisms that
represent individual matrices, then the composition of such morphisms by 
attaching directed arrows head-to-tail corresponds to the multiplication of these matrices. 
\vspace{3mm}

\noindent {\bf Remark.} It is also natural to represent $A$ to $End(V)$, endomorphisms of a 
given vector space $V$. The reader should note, however, that in this language the 
composition of morphisms of vector spaces proceeds in the opposite order from the 
matrix multiplication that we have preferred. In order to rectify this, one must speak of
representations of the opposite algebra $A^{op}$ to $End(V)$.
This point of view is useful in other contexts, but will not be pursued here.
\vspace{3mm}

\noindent {\bf Theorem.} Let $A$ be a (balanced) quantum algebra. Let $Rep:A \longrightarrow
M_{n}(k)$ be a representation of $A$ to the ring of $n \times n$ matrices over the ground 
ring $k$. Suppose there is an invertible  matrix $M$ in $M_{n}(k)$ such that the
automorphism $\tau:M_{n}(k) \longrightarrow M_{n}(k)$ given by $\tau(x) = MxM^{-1}$ gives 
$Rep(A)$ the structure of a (balanced) oriented quantum algebra such that 
$Rep(T(a)) = \tau(Rep(a))$ for all $a$ in $A$. 
Then the matrices $Rep(\rho)$, $Rep(\rho^{-1})$
and $M$ give data for a matrix model oriented link invariant of regular isotopy.
\vspace{3mm}

\noindent {\bf Proof.} The proof of this Theorem is contained in the discussion that 
precedes it. //
\vspace{3mm}

Some discussion of this Theorem is in order.  First of all, note that the evaluation of a
matrix model can be read directly from the pattern of matrices on the diagram. One simply
writes down the list of the diagram matrix entries corresponding to the cups, caps and 
crossings and then sums the product of the elements in this list over all possible values
of the repeated indices. The resulting evaluation is a function of the non-repeated indices.
\vspace{3mm}

\noindent On the other hand, we can view the matrices  that come from the representation of
an oriented quantum algebra as having the specific forms such as 

$$R^{ab}_{cd} = e'^{a}_{d}e^{b}.$$

By putting the expression of the invariant in terms of the signifier matrices $Rep(e)$,
$Rep(e')$, $Rep(E)$, $Rep(E')$ and $M$ and $M^{-1}$ we bring the expression of the invariant
close to the abstract expression in $Cat(A).$ In particular, there is a corresponding matrix
expression for any diagram that is obtained from the original abstract diagram by bead
sliding (just replace the signifiers by their corresponding matrices and write down the 
resulting sum of products of terms invloving them and the cups and caps). Thus every 
diagram obtained by sliding beads on the original abstract diagram has a matrix evaluation.
It is easy to see that these evaluations are invariant under bead sliding. (We leave the 
proof to the reader.) The consequence is that we can slide the beads first and then evaluate
the matrix model. (This pattern was first observed in \cite{GAUSS}.)
In particular if $K$ is a closed loop diagram and we slide all the algebra
to one segment where it takes the form $w$ and regularly homotop the flat diagram to the 
form $G^{n}$ then the matrix model yields the invariant 

$$INV(K) = trace((M^{2})^{n}Rep(w))$$

as the value of the matrix model. Here $trace$ denotes  standard matrix trace and the term
$M^{2}$ is the matrix evaluation of $G$ (Bead sliding shows that $GxG^{-1} = UD(x)$ for all
matrices $x$. The formula $G=M^{2}$ then follows from $U(x) = D(x) = MxM^{-1}.$) 
\vspace{3mm}

\noindent {\bf Remark on General Matrix Models.} Everything that we have said in this 
section generalizes to matrix models with two background matrices corresponding to the 
two automorphisms $U$ and $D$ in the general case of oriented quantum algebras. We have 
restricted ourselves to the balanced case only for ease of exposition. The general case 
follows the same lines. Since, in the genral case there are two automorphisms, there will
be two background matrices $M$ and $M'$ with $U(x) = MxM^{-1}$ and $D(x) = M'xM'^{-1}$ for
all matrices $x$. Specifically we will have

$$\overline{M}^{>} = M, \underline{M}_{>} = M^{-1}$$

\noindent and 

$$\overline{M}^{<} = M'^{-1}, \underline{M}_{<} = M'.$$

\noindent Otherwise, the model for a general oriented quantum algebra behaves in all 
respects like the model for a balanced algebra.
\vspace{3mm}

\subsection{An Example of a Matrix Oriented Balanced Quantum Algebra}

In this subsection we give a specific example of a matrix oriented balanced
quantum algebra. This algebra can be used to produce a sequence of
models of specializations of the Homfly polynomial \cite{KNOTS} \cite{JONES}, as
we shall see in the next section. \vspace{3mm}

\noindent Recall from the previous section the elementary matrices 
$E^{b}_{a}$ and the rule of multiplication $E_{a}^{b}E_{c}^{d} = \delta(b,c)E_{a}^{d}.$
\vspace{3mm}

\noindent We let 

$$z = q - q^{-1}.$$

\noindent In the equations below there is an implicit summation over the repeated
indices for an index set of the form $\{ 1,2,3,...,N\}$ for a natural number $N$.
Logical conditions on the indices are expressed within the formulas. The logical symbol
refers to indices at the same level in the expression to which it belongs.
For example,

$$E^{a}_{b} \otimes^{>} E^{b}_{a} = \Sigma_{a > b}E^{a}_{b} \otimes E^{b}_{a}.$$

\noindent We begin by defining $\rho$ , $\rho^{-1}$ and the automorphism $T.$
\vspace{3mm}

$$\rho = z E^{a}_{b} \otimes^{>} E^{b}_{a} +q E^{a}_{a} \otimes
E^{a}_{a} + E^{a}_{a} \otimes^{\neq} E^{b}_{b}$$

$$\rho^{-1} = -z E^{a}_{b} \otimes^{>} E^{b}_{a} +q^{-1} E^{a}_{a}
\otimes E^{a}_{a} + E^{a}_{a} \otimes^{\neq} E^{b}_{b}$$

$$T(E^{a}_{b}) = q^{a-b}E^{a}_{b}$$

\noindent Then

$$(1 \otimes T) \rho = zq^{a-b} E^{a}_{b} \otimes^{>} E^{b}_{a} +q
E^{a}_{a} \otimes E^{a}_{a} + E^{a}_{a} \otimes^{\neq} E^{b}_{b}$$

$$(T \otimes 1) \rho^{-1} = -zq^{b-a} E^{a}_{b} \otimes^{>}
E^{b}_{a} +q^{-1} E^{a}_{a} \otimes E^{a}_{a} + E^{a}_{a}
\otimes^{\neq} E^{b}_{b}.$$

Here is the calculation:

\begin{eqnarray*}
[(1 \otimes T)\rho][(T \otimes 1)\rho^{-1}] = 
-\Sigma_{a>b,a'>b'}z^{2}q^{a-b}q^{b'-a'}E^{a}_{b}E^{a'}_{b'} \otimes E^{b'}_{a'}E^{b}_{a} \\ 
+\Sigma_{a>b}zq^{a-b}q^{-1}E^{a}_{b}E^{a'}_{a'} \otimes E^{a'}_{a'}E^{b}_{a} 
+\Sigma_{a>b,a' \ne b'}zq^{a-b}E^{a}_{b}E^{a'}_{a'} \otimes E^{b'}_{b'}E^{b}_{a} \\ 
-\Sigma_{a'>b'}zqq^{b'-a'}E^{a}_{a}E^{a'}_{b'} \otimes E^{b'}_{a'}E^{a}_{a}
+\Sigma_{a,a'}E^{a}_{a}E^{a'}_{a'} \otimes E^{a'}_{a'}E^{a}_{a} \\ 
+\Sigma_{a' \ne b'}qE^{a}_{a}E^{a'}_{a'} \otimes E^{b'}_{b'}E^{a}_{a}
-\Sigma_{a \ne b,a'>b'}zq^{b'-a'}E^{a}_{a}E^{a'}_{b'} \otimes E^{b'}_{a'}E^{b}_{b} \\ 
+\Sigma_{a \ne b}q^{-1}E^{a}_{a}E^{a'}_{a'} \otimes E^{a'}_{a'}E^{b}_{b}
+\Sigma_{a \ne b,a' \ne b'}E^{a}_{a}E^{a'}_{a'} \otimes E^{b'}_{b'}E^{b}_{b} \\ 
\end{eqnarray*}

\begin{eqnarray*}
= -\Sigma_{a'>a>b'} z^{2}q^{a-b'}q^{a-a'}E^{a'}_{b'} \otimes E^{b'}_{a'}
+\Sigma_{a'>b'}z q^{a'-b'-1}E^{a'}_{b'} \otimes E^{b'}_{a'}  \\
-\Sigma_{a'>b'}z q^{b'-a'+1}E^{a'}_{b '} \otimes E^{b'}_{a'} 
+\Sigma_{a}E^{a}_{a} \otimes E^{a}_{a}
+\Sigma_{a \ne b} E^{a}_{a} \otimes E^{b}_{b}  \\
\end{eqnarray*}

\begin{eqnarray*}
= \Sigma_{a'>b'}z[-z \Sigma_{a'>a>b'} q^{a-b'}q^{a-a'} + q^{a'-b'-1}-q^{b'-a'+1}]
E^{a'}_{b'} \otimes E^{b'}_{a'} \\
+ \Sigma_{a} E^{a}_{a} \otimes E^{a}_{a}
+ \Sigma_{a \ne b} E^{a}_{a} \otimes E^{b}_{b} \\
\end{eqnarray*}

\noindent But

$$[-z \Sigma_{a'>a>b'} q^{a-b'}q^{a-a'} + q^{a'-b'-1}-q^{b'-a'+1}]
= q^{-a'-b'}[-z\Sigma_{a'>a>b'}q^{2a}+q^{2a'-1} -q^{2b'+1}]$$

\noindent and 
$$-z\Sigma_{a'>a>b'}q^{2a} = (q^{-1}-q)(q^{2b'+2} +...+q^{2a'-2}) = q^{2a'-1}-q^{2b'+1}$$

\noindent Thus

$$[(1 \otimes T)\rho][(T \otimes 1)\rho^{-1}] = \Sigma_{a}E^{a}_{a} \otimes E^{a}_{a}
+\Sigma_{a \ne b} E^{a}_{a} \otimes E^{b}_{b} = 1_{A \otimes A^{op}}.$$

This verifies that the algebra generated by elementary matrices in conjunction with 
$\rho$, $\rho^{-1}$ and $T$ form a balanced oriented quantum algebra. In order to associate
a link invariant to this algebra we can construct a matrix model by taking

$$M_{ij} = q^{i}\delta_{ij}.$$

\noindent We leave it as an exercise for the reader to verify that the resulting invariant
is an unnormalized version of the Homfly polynomial. This exercise is clarified by the 
remarks on state models in the next section.
\vspace{3mm}

\section{State Sums}

There are many versions of the general notion of a state summation model for a 
link invariant. Given a link diagramm $K$, a {\em combinatorial structure
associated with $K$} is a graph with decorations (possibly algebraic) that is 
obtained from $K$ by some well-defined process of labelling and replacement.
The exact details of what a combinatorial structure can be are left open, as there are
many possibilities. One well known way to obtain combinatorial structures associated to a
link diagram $L$ is to replace each crossing of $L$ with either a {\em smoothing} of that 
crossing or a {\em flattening} of that crossing. In smoothing a crossing we replace the
connections at the crossing so that it fits in the plane without any arcs crossing over or 
under one another. In flattening a crossing, the crossing is replaced by a four-valent
vertex in the plane. See Figure 20 for illustrations of smoothing and flattening.  
\vspace{3mm}  

By a  {\em combinatorial state sum} I mean that for each link diagram
$K$ there will be associated a set of combinatorial structures $S(K)$,
called the {\em states} of $K$, and a 
functional $<K|S>$ that associates to each state of $K$ and state $S$ an element $<K|S>$ 
in a commutative ring $k$. Then we define the {\em state sum for $K$} to be the summation 

$$ Z_{K} = \Sigma_{S \in S(K)} <K|S> <S>.$$

\noindent 
where $<S>$ is a specifically given state evaluation in $k$.
It is intended that $Z_{K}$ be a regular isotopy invariant of the link $K$. 
\vspace{3mm}

In this section we shall consider state sums of the following special form.
There is given an index set $I=\{ 1,2,3,...,n \}$. Each crossing of the link diagram 
can be replaced by either oriented parallel arcs (an oriented smoothing)
with a sign of either equality or
inequality between them (this constitutes three possibilities, since the inequality can
be oriented in two ways between the two lines), or by crossed arcs (a flattening)
with a sign of inequality between them. In Figure 20 we have illustrated these local
replacements in the form of symbolic summations (we use the diagram and its evaluation
interchangeably)

$$K_{+} = ASK_{=} + BSK_{<} +CSK_{>} + DFK_{\ne}$$
$$K_{-} = A'SK_{=} + B'SK_{<} +C'SK_{>} + D'FK_{\ne}$$

\noindent where $K_{+}$ and $K_{-}$ denote the link with a specific crossing that is either
positive or negative, $SK_{=}$ denotes the result of smoothing with left top line equal
to the left top right line,
$SK_{<}$ denotes the result of smoothing with left top line less than the 
left top right line, $SK_{>}$ denotes the result of smoothing with left top line greater 
than the 
left top right line, $FK_{\ne}$ denotes the result of flattening with left top line unequal 
to the
left top right line. These equations express the state summation symbolically. We
expand in this way on eacn crossing until a formal sum of products is reached. The diagrams
that are implicated in these products have a medley of stipulations of equality and
inequality inscribed upon them. Those that are impossible (e.g. a line is asked to be less
than itself) are given value zero. A diagram $D$ with a possible assignments of equality and 
inequality is evaluated as $<D>$ as in the general state sum description of the previous 
paragraph. This recursive description of the state summation is identical to the description

$$ Z_{K} = \Sigma_{S \in S(K)} <K|S> <S>$$

\noindent where the states $S$ are the diagrams $D$ obtained by smoothing and 
flattening, and the evaluations $<K|S>$ denote the product of labels $A$,$A'$, ...,$D$,$D'$
that are implicit in the recursive expansion. That is, one should think of that state 
diagram $D$ as decorated not only with signs of equality and inequality, but also with
the alphabetic labels that are implicit in the recursive expansion. Then $<K|S>$ is the
product of the labels that is the coefficient of the corresponding diagram in the 
recursive expansion. 
\vspace{3mm}

We now turn to the specific definition of $<D>$ for a state diagram $D$. This diagram is 
labelled with equalities and inequalities that make it possible to create actual labellings
of its curves (one numerical label per curve) from the index set 
$I=\{ 1,2,3,...,n \}.$
Let $\Lambda$ denote an actual labelling of $D$ and $D(\Lambda)$ that particular labelling
of $D$. Let 

$$<D> = \Sigma_{\Lambda} <D(\Lambda)> = 
\Sigma_{\Lambda} q^{\Sigma_{i \in \Lambda} iRot_{D}(i)}$$

\noindent where $\Lambda$ runs over all admissible labellings of $D$ from the index set $I$.
Here $q$ is a commuting algebraic variable, and $Rot_{D}(i)$ denotes the Whitney degree
in the plane of the component of $D$ that is labelled by the index $i$. (We think of 
$\Lambda$ as a list of indices (with multiplicities) such that each index is attached to 
a specific component of $D$.) 
\vspace{3mm}

This specification completely defines the state summation.
That is, the state summation is now well-defined on link diagrams. It is not in general 
an invariant of regular isotopy for links.  We will now show that this specification 
corresponds to a specific matrix model and that a special case gives specializations of 
the Homfly polynomial.
\vspace{3mm}

The Homfly specialization has vertex weights given by the expansion

$$K_{+} = qSK_{=} + (q-q^{-1})SK_{<}  + FK_{\ne}$$
$$K_{-} = q^{-1}SK_{=} + (q^{-1} - q)SK_{>} +FK_{\ne}.$$

\noindent Note that 

$$K_{+} - K_{-} = (q-q^{-1})[SK_{=} + SK_{<} + SK_{>}] = (q-q^{-1})SK$$

\noindent where $SK$ denotes the diagram obtained by smoothing the crossing in
question. This is the source of the skein relation for the Homfly polynomial.
\vspace{3mm}
   
\noindent The matrix model arises by interpreting the expansion formulas for 
$K_{+}$ and $K_{-}$
as definitions for the braiding matrices in the matrix model. Translating these braiding 
matrices to the algebraic form (by composing with the appropriate permutation) yields the
$\rho$ and $\rho^{-1}$ of the previous section. We omit these details. The upshot is that
the present state sum can be used to prove that the matrix model in the last section does
give the Homfly polynomial specializations.
\vspace{3mm}

In Figure 21 we show how to 
interpret the smoothed and flattened local states as matrices, using the convention that 
$\delta[a,b,c,d]$ is equal to 1 when $a=b=c=d$ and is equal to 0 otherwise,  
$\delta[a,b]$ is equal to 1 when $a=b$ and is equal to 0 otherwise and $[P]$ is  equal to 
1 when the proposition $P$ is true and is equal to 0 otherwise. With these interpretations 
for the braiding matrices, it is not hard to see that this state model is identical to the 
link invariant derived from the oriented quantum algebra of the previous section.
This completes our discussion of these relationships among state models, matrix models
and oriented quantum algebras.
\vspace{3mm}

 {\tt    \setlength{\unitlength}{0.92pt}
\begin{picture}(428,232)
\thinlines    \put(338,47){D'}
              \put(262,48){C'}
              \put(191,49){B'}
              \put(319,174){D}
              \put(254,175){C}
              \put(181,175){B}
              \put(347,191){\framebox(38,30){$\ne$}}
              \put(358,69){\framebox(38,30){$\ne$}}
              \put(16,137){\vector(3,4){63}}
              \put(86,142){\vector(-1,1){33}}
              \put(43,185){\vector(-1,1){33}}
              \put(349,11){\vector(3,4){63}}
              \put(417,15){\vector(-1,1){76}}
              \put(123,140){\vector(0,1){77}}
              \put(155,140){\vector(0,1){78}}
              \put(225,139){\vector(0,1){78}}
              \put(194,140){\vector(0,1){77}}
              \put(293,139){\vector(0,1){78}}
              \put(265,139){\vector(0,1){77}}
              \put(90,22){\vector(-1,1){76}}
              \put(18,22){\vector(1,1){30}}
              \put(59,63){\vector(1,1){34}}
              \put(281,14){\vector(0,1){77}}
              \put(309,14){\vector(0,1){78}}
              \put(210,14){\vector(0,1){77}}
              \put(238,13){\vector(0,1){78}}
              \put(162,13){\vector(0,1){78}}
              \put(133,13){\vector(0,1){77}}
              \put(404,140){\vector(-1,1){76}}
              \put(336,136){\vector(3,4){63}}
              \put(199,170){\framebox(21,20){<}}
              \put(270,168){\framebox(19,21){>}}
              \put(138,41){\framebox(20,19){=}}
              \put(130,167){\framebox(20,19){=}}
              \put(214,43){\framebox(21,20){<}}
              \put(285,42){\framebox(19,21){>}}
              \put(80,165){\makebox(38,23){= A}}
              \put(160,166){\makebox(25,23){+}}
              \put(231,165){\makebox(21,23){+}}
              \put(297,164){\makebox(27,24){+}}
              \put(103,42){\makebox(27,21){= A'}}
              \put(166,39){\makebox(26,25){+}}
              \put(243,39){\makebox(23,22){+}}
              \put(313,37){\makebox(30,26){+}}
\end{picture}}

\noindent
{\bf Figure 20 - Local States}
\vspace{3mm}

{\tt    \setlength{\unitlength}{0.92pt}
\begin{picture}(381,269)
\thinlines    \put(224,80){\framebox(38,30){$\ne$}}
              \put(13,165){\vector(0,1){77}}
              \put(45,165){\vector(0,1){78}}
              \put(248,166){\vector(0,1){78}}
              \put(217,167){\vector(0,1){77}}
              \put(44,35){\vector(0,1){78}}
              \put(16,35){\vector(0,1){77}}
              \put(281,29){\vector(-1,1){76}}
              \put(213,25){\vector(3,4){63}}
              \put(222,197){\framebox(21,20){<}}
              \put(21,64){\framebox(19,21){>}}
              \put(20,192){\framebox(20,19){=}}
              \put(10,250){a}
              \put(43,249){b}
              \put(11,153){c}
              \put(44,152){d}
              \put(54,198){$= \delta[a,b,c,d]$}
              \put(215,251){a}
              \put(246,250){b}
              \put(214,156){c}
              \put(245,154){d}
              \put(259,202){$= \delta[a,c] \delta[b,d] [a < b]$}
              \put(14,118){a}
              \put(41,119){b}
              \put(13,23){c}
              \put(45,22){d}
              \put(55,71){$= \delta[a,c] \delta[b,d] [a > b]$}
              \put(203,114){a}
              \put(274,115){b}
              \put(210,13){c}
              \put(278,17){d}
              \put(287,68){$ \delta[a,d] \delta[b,c] [a \ne b]$}
\end{picture}}

\noindent
{\bf Figure 21 - Local Matrices}
\vspace{3mm}

\section{Quasitriangular Ribbon Hopf Algebras} 

In this section we show how to recover the results of Reshetikhin 
and Turaev \cite{RIBBON,RESH}  that
associate an oriented link invariant to each representation of a
quasitriangular ribbon Hopf algebra. In our language, each such
representation gives rise to a $D=1$  (standard)
oriented quantum algebra and a corresponding matrix model. 
\vspace{3mm}

Any quasitriangular Hopf algebra is an example of an
(unoriented) quantum algebra. In a quantum algebra with 
antipode $s$ we have

$$(s \otimes s) \rho = \rho$$

\noindent and 

$$(s \otimes 1) \rho = (1 \otimes s^{-1}) \rho = \rho^{-1}.$$

\noindent Then

$$(se)f \otimes e'f' = 1 \otimes 1$$

$$(sf)s^{2}(e) \otimes e'f' = 1 \otimes 1$$

$$(sf)e \otimes s^{-2}(e')f' = 1 \otimes 1$$

$$[(s \otimes 1) \rho][(1 \otimes s^{-2}) \rho] = 1_{A} \otimes
1_{A^{op}}$$

$$[\rho^{-1}][(1 \otimes s^{-2}) \rho] = 1_{A} \otimes 1_{A^{op}}.$$

\noindent Similarly,

$$[(1 \otimes s^{-2}) \rho][\rho^{-1}] = 1_{A} \otimes 1_{A^{op}}.$$

\noindent This shows that any quantum algebra $A$ is a standard
oriented quantum algebra with $t=s^{-2}$.  In particular, This
applies to quasitriangular ribbon Hopf algebras. We
now obtain in this way a new proof of the reuslts of Reshetikhin
and Turaev on the existence of link invariants from representations of
such algebras. 
\vspace{3mm}

In order to define a link invariant from a representation, 
we will define cups and caps as morphisms of vector spaces 
so that they represent the automorphisms $t=U$ and $1=D.$ Once
these cups and caps are defined, the evaluation of tangles is accomplished by regarding 
the image of each tangle under the functor $F$ as a morphism of vector spaces. In 
particular, a closed link diagram is a linear map from $k$ to $k$ and the value of its
invariant is the value of this map on the unit $1$ in $k$. This method of evaluation
is neccessary for the direct comparison of our invariant with the Reshetikhin Turaev 
invariant. It is also compatible with the other methods in this paper (matrices and 
bead sliding).
\vspace{3mm}

In a quasitrangular ribbon Hopf algebra $A$ the square of the 
antipode is represented by a grouplike element $G$ \cite{DRIN} 
so that $s^{2}(a) = G^{-1}aG$ and 
$$t(a) = s^{-2}(a) = GaG^{-1}$$ 
\noindent for all $a$ in $A$.  The existence of this grouplike
element allows us, given a finite dimensional representation
of $A$, to represent the cups and the caps so that they correctly
represent the automorphisms $t$ and the identity. The definitions
are given below.
\vspace{3mm}

\noindent {\bf Remark.} A quantum algebra with anti-automorphism
$s$ whose square $s^{2}$ is represented by conjugation by an element $G$ is called 
(by us) a {\em twist quantum algebra} \cite{RKinv}. The results of this section can be 
stated in full generality for twist quantum algebras.
\vspace{3mm} 

Since we are in a
representation of the algebra $A$, we can assume that each 
element of $A$ corresponds to an endomorphism of a vector
space $V$. Let $\beta$ run over a basis $\cal B$ for
$V$. Let $\beta^{*}$ run over the dual basis for $V^{*}$. 
Thus $\beta^{*}(\beta') = \delta(\beta,\beta').$ 
\vspace{3mm}

Note that for any $v$ in $V$ that

$$v = \Sigma_{\beta} \beta^{*}(v)\beta.$$
\vspace{3mm}

\noindent We define
$v^{*}$ for any $v$ in $V$ by the equation 
$v^{*}(\beta) = \beta^{*}(v)$  where $\beta$ is in the basis $\cal B.$

\vspace{3mm}

The transpose of an element $a$ in $A$ will be denoted
by $a^{t}.$ Then $a^{t}$ corresponds to an element in the 
endomorphisms of the dual space $V^{*}.$  
Note also that if $a$ is an endomorphism of $V$, then $a^{t}$ is an endomorphism
of $V^{*}$ defined by the formula $$a^{t}f(v) = f(a(v))$$ 

\noindent for any $f$ in $V^{*}$.
By our usual conventions
the transpose of the morphism for an element $a$ is diagrammed 
by drawing a down line that is still labelled $a$.  
\vspace{3mm}

First we define $CapRight$ and $CupRight.$ It is understood that
these formulas occur in the representation with $x \otimes y^{*}$
representing an element of $V\otimes V^{*}$. 
Each summation runs over a basis for $V$.
\vspace{3mm} 
 
$$CapRight(x \otimes y^{*}) = y^{*}(G^{-1}x)$$

\noindent for any $x$ and $y$ in $V.$ 

$$CupRight(1) = \Sigma_{\beta}(\beta^{*} \otimes G\beta).$$

\noindent where $\beta$ runs over a basis for $V$.

\noindent Then, letting $a$ be any endomorphism of $V$ and 
$x$ an element of $V$, we let $t(a)$ be the 
endomorphism of $V$ obtained from sandwiching $a$ in the middle
of the CapRight and CupRight as in Figure 10. We have

$$t(a)(x) = \Sigma_{\beta} CapRight(x \otimes a^{t}\beta^{*})G\beta 
= \Sigma_{\beta} a^{t}\beta^{*}(G^{-1}x)G\beta
= \Sigma_{\beta} \beta^{*}(aG^{-1}x)G\beta$$

\noindent Now

$$aG^{-1}x = \Sigma_{\beta} \beta^{*}(aG^{-1}x)\beta$$

\noindent thus

$$t(a)x = GaG^{-1}x.$$

\noindent hence $$t(a) = GaG^{-1}.$$ 

\noindent
$CapLeft$ and $CupLeft$ are defined by evaluation and coevaluation
without the use of $G$ so that $D=1$. \vspace{3mm}

These constructions show that any representation of a 
quasitriangular ribbon Hopf algebra $A$ or an oriented twist quantum algebra
\cite{RKinv} gives rise to an invariant
of regular isotopy of knots and links. The specification of the 
cups and caps gives a matrix model for this invariant in any 
basis for the representation space. In fact, the specification of cups and caps that 
we have used matches (up to reversals of arrow conventions)
the cups and caps of Reshetikhin and Turaev in \cite{RIBBON}.
The choice of cups and caps determines the automorphisms of the corresponding 
oriented quantum algebra. The oriented quantum algebra underlying the
Reshetikhin Turaev invariants for a representation of a ribbon Hopf algebra is identical to 
the algebra that we have described in this section. From this it follows that these methods
reproduce the Reshetikhin Turaev invariants.
 
 \end{document}